\documentclass[francais,psamsfonts]{smfart}
\usepackage[english,francais]{babel}
\usepackage{smfthm,bull,epsfig,smfenum}

\numberwithin{equation}{section}                

\newcommand{\sst}{\scriptscriptstyle}               

\newcommand{\C}{{\mathbb C}}  
\newcommand{\Z}{{\mathbb Z}}
\newcommand{\N}{{\mathbb N}}     
\newcommand{\R}{{\mathbb R}}          

\newcommand{\id}{\mathop{\rm id}\nolimits}                                    
\newcommand{\End}{\mathop{\rm End}\nolimits}                    

\newcommand{\hdx}{U_h(D_x)}

\newcommand{\ad}{{{\rm ad}_+}}

\newcommand{\e}[1]{E_{\beta_{#1}}}
\newcommand{\f}[1]{F_{\beta_{#1}}}
\renewcommand{\k}[1]{K_{\beta_{#1}}}

\def\mybaselines{\baselineskip10pt%
    \lineskip3pt \lineskiplimit3pt}

\catcode`@=11
\def\smallmatrixp#1{\null\,\vcenter{\mybaselines\m@th
    \ialign{\hfil$\scriptscriptstyle##$\hfil&&\!\hfil$\scriptscriptstyle##$
    \hfil\crcr\mathstrut\crcr\noalign{\kern-\baselineskip}
    #1\crcr\mathstrut\crcr\noalign{\kern-\baselineskip}}}\kern-2pt}
\def\psmallmatrixp#1{\left(\smallmatrixp{#1}\right)}
\catcode`@=12

\catcode`@=11
\def\smallmatrixg#1{\null\,\vcenter{\mybaselines\m@th
    \ialign{\hfil$\scriptscriptstyle##$\hfil&&\space\hfil$\scriptscriptstyle##$
    \hfil\crcr\mathstrut\crcr\noalign{\kern-\baselineskip}
    #1\crcr\mathstrut\crcr\noalign{\kern-\baselineskip}}}\kern-2pt}
\def\psmallmatrixg#1{\left(\smallmatrixg{#1}\right)}
\catcode`@=12

\newcommand{\dessin}[1]{\vcenter{\hbox{\epsfbox{#1}}}}

\newcounter{fig}

\author{Henrik Thys}

\address{D\'epartement de Math\'ematiques - IRMA\\
         7, rue Ren\'e Descartes \\
         67084 Strasbourg C\'edex }
         
\email{thys@math.u-strasbg.fr}

\title[$U_h(D(2,1,x))$, $R$-matrice universelle et invariant d'entrelacs]%
   {$R$-matrice universelle pour $U_h(D(2,1,x))$ %
       et invariant d'entrelacs associ\'e}
       
\alttitle{Universal $R$-matrix for $U_h(D(2,1,x))$ %
          and link invariant arising from it}

\begin{document}

\frontmatter

 
\begin{abstract}
En utilisant
la m\'ethode du double quantique, nous construisons une $R$-matrice universelle 
pour la quantification de la superalg\`ebre de Lie $D(2,1,x)$. Nous utilisons
ce r\'esultat pour construire un invariant d'entrelacs et nous montrons qu'il
est \'egal \`a une sp\'ecialisation du polyn\^ome de Dubrovnik introduit par 
Kauffman.
\end{abstract}

\begin{altabstract}
Using the quantum double method, we construct a universal $R$-matrix for
the quantization of the Lie superalgebra $D(2,1,x)$. We use this
result to construct a link invariant and show it coincides with a 
specialization of Kauffman's Dubrovnik polynomial.
\end{altabstract}



\altkeywords{quantum supergroup, universal $R$-matrix, quantum double,
             knots invariant, Lie superalgebra}
             
\maketitle                       


\mainmatter


\section*{Introduction}
\addcontentsline{toc}{section}{\numberline{}Introduction}

Les groupes quantiques introduits autour de 1983-85 par Drinfeld et Jimbo
sont des d\'eformations \`a un param\`etre des alg\`ebres enveloppantes des 
alg\`ebres de Lie semisimples complexes. Techniquement, ces ``quantifications" 
sont des alg\`ebres de Hopf munies de ``$R$-matrices universelles", 
c'est-\`a-dire d'\'el\'ements qui sont responsables de l'existence de solutions 
de la fameuse \'equation de Yang-Baxter et donc de repr\'esentations des groupes 
de tresses.

Dans les ann\'ees 1970 Victor Kac \cite{kac} a \'etudi\'e une g\'en\'eralisation 
naturelle 
des alg\`ebres de Lie semisimples, \`a savoir les superalg\`ebres de Lie. Dans 
la classification qu'il en donne, il y a celles que l'on peut appeler des 
superalg\`ebres ``classiques" comme $sl(n | m)$ ou $osp(n | m)$ et il y a des 
superalg\`ebres exceptionnelles. Parmi ces derni\`eres, il y a une famille, 
not\'ee 
$D(2,1,x)$ dans~\cite{kac}, et que nous noterons $D_x$, d\'ependant d'un 
param\`etre continu~$x$.
La superalg\`ebre de Lie $D_x$ joue un r\^ole particulier en physique o\`u
elle fournit la seule th\'eorie topologique quantique des champs (TQFT) de 
Chern-Simons pour laquelle la th\'eorie conforme des champs (CFT) de dimension 
deux correspondante a une supersym\'etrie~$N = 4$. Elle joue aussi un r\^ole
particulier dans les travaux r\'ecents de Pierre Vogel \cite{vogel}
et de Jens Lieberum \cite{jens} sur les invariants de Vassiliev.

Apr\`es les alg\`ebres de Lie semisimples, les superalg\`ebres de Lie
ont \'egale\-ment \'et\'e quantifi\'ees (par Gould {\em et al.}, Leites 
{\em et al.}, Scheunert, etc., {\em cf.} par exemple~\cite{CK,FLV,S93,Y96}). Les quantifications obtenues sont
des superalg\`ebres de Hopf munies de bases de type Poincar\'e-Birkhoff-Witt.

Pour ce qui est de l'existence d'une $R$-matrice universelle pour les 
supergroupes quantiques, elle a \'et\'e \'etablie pour les quantifications
de toutes les superalg\`ebres de Lie classifi\'ees par Kac, \`a l'exception 
pr\'ecis\'ement de~$D_x$. Dans cette partie, nous comblons cette lacune
en construisant explicitement une $R$-matrice universelle pour la 
quantification~$\hdx$.

La m\'ethode utilis\'ee est celle du double quantique introduite par
Drinfeld, m\'ethode dont se sont servis Rosso~\cite{rosso89},
Kirillov-Reshetikhin~\cite{kirillov-reshet90} et 
Levendorsky-Soibelman~\cite{leven-soib90} pour les groupes quantiques. Cette m\'ethode 
s'\'etend au cas $\Z/2\Z$-gradu\'e. Nous d\'efinissons des analogues 
des vecteurs de racines positives et n\'egatives pour $\hdx$. Nous 
d\'efinissons \'egalement
l'analogue de la sous-alg\`ebre positive $U_+$ (resp. sous-alg\`ebre n\'egative
$U_-$) engendr\'ee par les vecteurs de racines positives (resp. n\'egatives), et
nous construisons un accouplement de Hopf entre $U_+$ et~$U_-$. Ensuite, nous 
calculons les  relations de commutation entre les vecteurs de racines ainsi que 
leur coproduit. Ceci nous permet de construire des bases de $U_+$ et~$U_-$, 
duales pour l'accouplement de Hopf. Nous en d\'eduisons une $R$-matrice 
universelle pour~$\hdx$.

Une $R$-matrice universelle sur une quantification formelle munit la cat\'egorie 
des modules topologiques d'une structure de cat\'egorie ruban\'ee au sens de
Turaev, {\em cf.} \cite{K,Z95}. Dans le cas qui nous int\'eresse, ceci fournit 
pour chaque 
$\hdx$-module un invariant d'isotopie d'entrelacs parall\'elis\'es et orient\'es. 
En prenant un supermodule de dimension six, nous obtenons un invariant 
${\mathcal I}$ d'entrelacs
parall\'elis\'es non orient\'es, \`a valeurs dans $\Z[q,q^{-1}]$, et nous montrons
qu'il est \'egal \`a une sp\'ecialisation du polyn\^ome de Dubrovnik introduit par
Kauffman.

Le plan est le suivant. Les \S\S\ref{rappelshopf} et \ref{hdx} sont 
consacr\'es \`a des rappels sur les superalg\`ebres de Hopf et sur 
$\hdx$. Au \S\ref{resultats} nous \'enon\c{c}ons les r\'esultats principaux
(th\'eor\`emes~\ref{theoR} et~\ref{theo2}). Le \S\ref{preliminaire}
est consacr\'e \`a la construction d'un accouplement de Hopf et de ses 
propri\'et\'es
et le \S\ref{D} \`a la construction d'un double quantique ${\mathcal D}$. Au 
\S\ref{relations_D} nous \'etablissons des relations v\'erifi\'ees dans ${\mathcal D}$ et
nous donnons la d\'emonstration du th\'eor\`eme  
\ref{theoR} au \S\ref{demo_dx}. Au \S\ref{par:module}, nous d\'efinissons un
$\hdx$-module $M$ de rang six et nous calculons le tressage correspondant 
\`a l'aide du
th\'eor\`eme~\ref{theoR}. Au \S\ref{catrub}, nous \'enon\c{c}ons quelques
propri\'et\'es de la cat\'egorie ruban\'ee associ\'ee \`a $M$, et nous
terminons par la d\'emonstration du th\'eor\`eme~\ref{theo2} au \S\ref{invariant}.

Cet article est tir\'e du troisi\`eme chapitre de ma th\`ese \cite{moi}. Je
tiens \`a remercier J. Alev qui m'a encourag\'e \`a construire un invariant
de n\oe ud \`a partir de la $R$-matrice du~\S\ref{par:module}, et C. Kassel pour son
aide.


\section{\'Enonc\'e des th\'eor\`emes principaux}
\label{theoprinc}

Dans ce paragraphe, nous commen\c{c}ons par quelques rappels sur les 
superalg\`ebres de Hopf, l'\'equation de Yang-Baxter gradu\'ee et les 
$R$-matrices universelles. Nous donnons ensuite la d\'efinition
de la superalg\`ebre de Hopf $\hdx$ qui est la quantification
de la superalg\`ebre de Lie $D_x$. Nous terminerons par l'\'enonc\'e des 
th\'eor\`emes principaux.

Le contenu du \S\ref{rappelshopf} se trouve dans de nombreux articles. On pourra
notamment consulter \cite{CK,FLV,S93,Y94a,Y96,Z92}. 

\subsection{Superalg\`ebres de Hopf tress\'ees et double quantique}
\label{rappelshopf}

On note $\C$ le corps des nombres complexes. Un {\em superespace vectoriel} $V$ 
est un $\C$-espace vectoriel muni d'une graduation par
$\Z/2\Z$, {\em i.e.} d'une somme directe de deux espaces vectoriels
$V=V_0\oplus V_1,$
o\`u $V_0$ est la partie {\em paire} de $V$ et $V_1$ la partie {\em impaire}. 
Les \'el\'ements de
$V_0$ (resp. de $V_1$) sont dits {\em homog\`enes pairs} (resp. {\em homog\`enes 
impairs}). Si $v\in V_0$, on pose $|v|=0$ et si $v\in V_1$, on pose 
$|v|=1$, et on appelle {\em degr\'e} ces quantit\'es. Le corps $\C$ est un
superespace dont la partie impaire est r\'eduite \`a $0$.
Le produit tensoriel de deux superespaces $V$ et $W$ est le superespace
$V\otimes W=(V\otimes W)_0\oplus (V\otimes W)_1,$
o\`u
$$(V\otimes W)_0=(V_0\otimes V_0)\oplus (V_1\otimes V_1),\quad
(V\otimes W)_1=(V_0\otimes V_1)\oplus (V_1\otimes V_0).$$
\'Etant donn\'e deux superespaces $V$ et $W$, un {\em morphisme de 
superespaces}
$f:V\rightarrow W$
est une application lin\'eaire telle que 
$f(V_i)\subset W_i.$
Dans toute la suite, le morphisme identit\'e d'un superespace $V$ sera not\'e
$\id_V$.
La {\em volte} de deux superespaces $V,W$ est le morphisme de superespaces
$\tau :V\otimes W\longrightarrow W\otimes V$
d\'efini sur des \'el\'ements homog\`enes $v\in V,w\in W$ par
$$
\tau (v\otimes w)=(-1)^{|v||w|}w\otimes v.
$$
(La volte sera not\'ee $\tau$ pour tous les couples de superespaces.)

Une {\em superalg\`ebre} est un triplet $(A,\mu ,\eta )$  o\`u $A$ est un 
superespace, $\mu :A\otimes A\rightarrow A$ et $\eta :\C\rightarrow A$ des 
morphismes de superalg\`ebres
tels que $\mu (\mu\otimes\id_A)=\mu (\id_A\otimes\mu)$ et
 $\mu (\eta\otimes\id_A)=\mu (\id_A\otimes\eta )=\id_A$.
On notera $\mu (a\otimes a')=aa'$ pour $a,a'\in A$. \'Etant donn\'e deux 
alg\`ebres $(A,\mu_A ,\eta _A)$ et $(B,\mu_B ,\eta_B )$, un morphisme de
superalg\`ebres $f:A\rightarrow B$ est un morphisme de superespaces tel que
$f(aa')=f(a)f(a')$ pour tous $a,a'\in A$.
Le produit tensoriel de deux superalg\`ebres $(A,\mu_A ,\eta _A)$ et
$(B,\mu_B ,\eta_B )$ est une superalg\`ebre $(A\otimes B,\mu_{A\otimes B},
\eta_{A\otimes B})$ o\`u
$\mu_{A\otimes B}=(\mu_A\otimes\mu_B)\circ (\id_A\otimes\tau\otimes\id_B)
\text{ et } \eta_{A\otimes B}=\eta_A\otimes\eta_B.$
On notera que pour des \'el\'ements homog\`enes $a,a'\in A$ et $b,b'\in B$, on a
\begin{equation}\label{produit}
(a\otimes b)(a'\otimes b')=(-1)^{|b||a'|}aa'\otimes bb'.
\end{equation}
Une {\em superalg\`ebre de Hopf} $(A,\mu ,\eta ,\Delta, \varepsilon ,S)$ est  
la donn\'ee 
d'une superalg\`ebre $(A,\mu ,\eta )$, de morphismes de superalg\`ebres 
$\Delta :A\rightarrow A\otimes A$ 
({\em le coproduit}) et
$\varepsilon :A\rightarrow \C$ ({\em la co\"unit\'e}), et d'un morphisme de 
superespaces $S:A\rightarrow A$ ({\em l'antipode})  tels que
\begin{enumerate}
\item  $(\Delta\otimes\id_A)\Delta =(\id_A\otimes\Delta )\Delta$,
\item  $(\varepsilon\otimes\id_A)\Delta =(\id_A\otimes\varepsilon )\Delta=
\id_A$,
\item  $\mu (\id_A\otimes S)\Delta =\mu (S\otimes\id_A)\Delta=
      \eta\varepsilon $.
\end{enumerate}
Rappelons que l'antipode est un anti-morphisme de superalg\`ebre, {\em
i.e.}
$$S(aa')=(-1)^{|a||a'|}S(a')S(a),\quad a,a'\in A\text{ homog\`enes}.$$

Une superalg\`ebre de Hopf $A$ est dite {\em tress\'ee} s'il existe un 
\'el\'ement inversible $R=\sum a_i\otimes b_i\in (A\otimes A)_0$ tel que
\begin{gather*}
R\Delta (a)=(\tau\circ\Delta )(a)R,\quad\forall\; a\in A,\\
(\Delta\otimes\id_A)R=R_{13}R_{23},\quad
(\id_A\otimes \Delta )R=R_{13}R_{12},
\end{gather*}
o\`u
$$R_{12}=\sum_ia_i\otimes b_i\otimes 1,\quad R_{13}=\sum_ia_i\otimes 1\otimes 
b_i,
\quad R_{23}=\sum_i 1\otimes a_i\otimes b_i\in A\otimes A\otimes A.$$
L'\'el\'ement $R$ est appel\'e {\em $R$-matrice universelle} de $A$. Il v\'erifie
l'\'equation de Yang-Baxter gradu\'ee 
$$R_{12}R_{13}R_{23}=R_{23}R_{13}R_{12}.$$

Nous rappelons maintenant la construction du {\em double quantique} de
Drinfeld. On pourra consulter \cite{drinfeld85,KRT,rosso89bis,ZGB} pour les
d\'etails.
Soient $A$ et $B$ deux superalg\`ebres de Hopf avec antipode inversible. Un
{\em accouplement de Hopf} est un morphisme de superespaces 
\mbox{$\varphi :B\otimes A\rightarrow \C$} tel que
\begin{equation}\label{accoup2}
\begin{gathered}
\varphi (bb'\otimes a)=\sum_{(a)}(-1)^{|b'||a_{(1)}|}\,\varphi 
(b\otimes a_{(1)})\varphi (b'\otimes a_{(2)}),\\
\varphi (b\otimes aa')=\sum_{(b)}\varphi 
(b_{(1)}\otimes a')\varphi (b_{(2)}\otimes a),\\
\varphi (1\otimes a)=\varepsilon (a)
\text{ et } \varphi (b\otimes 1)=\varepsilon (b),\\
\varphi (b\otimes S(a))=\varphi (
S^{-1}(b)\otimes a).
\end{gathered}
\end{equation}
pour tous les \'el\'ements homog\`enes  $a,a'\in A$ et $b,b'\in B$, et o\`u
$$\Delta (a)=\sum_{(a)} a_{(1)}\otimes a_{(2)},$$
suivant la notation de Sweedler.
On construit \`a partir de $A$ et $B$ une superalg\`ebre 
de Hopf ${\mathcal D}(A,B,\varphi)$ de la mani\`ere suivante :
\begin{enumerate}
\item ${\mathcal D}(A,B,\varphi)=A\otimes B$ comme superespace.
\item Le coproduit de ${\mathcal D}(A,B,\varphi)$ est donn\'e par 
$(\id_{A}\otimes\,\tau\otimes\id_{B})(\Delta_{A}\otimes
\Delta_{B})$.
\item La co\"unit\'e est le produit tensoriel des co\"unit\'es de $A$ et 
$B$.
\item L'unit\'e est le produit tensoriel des unit\'es de $A$ et $B$.
\item Le produit est d\'efini par les deux formules suivantes:
\begin{gather*}
(a\otimes 1)(1\otimes b)=a\otimes b,\\
(1\otimes b)(a\otimes 1)=
\sum_{(a),(b)}(-1)^\xi\,\varphi (S(b_{(1)})\otimes a_{(1)})\varphi 
(b_{(3)}\otimes a_{(3)})a_{(2)}\otimes b_{(2)},
\end{gather*}
o\`u $\xi =|a_{(1)}||b_{(2)}|+|a_{(2)}||b_{(2)}|+|a_{(1)}||b_{(3)}|
+|a_{(2)}||b_{(3)}|$ et $a\in A$, $b\in B$ sont homog\`enes.
\item L'application $a\mapsto a\otimes 1$ (resp. $b\mapsto 1\otimes b$)
de $A$ dans ${\mathcal D}(A,B,\varphi)$ (resp. de $B$ dans 
${\mathcal D}(A,B,\varphi)$) 
est un morphisme de superalg\`ebres de Hopf injectif.
\end{enumerate}
Ces deux derniers points nous permettent d'identifier $a$ \`a $a\otimes  1$ 
pour $a\in A$ et
$b$ \`a $1\otimes b$ pour $b\in B$, et ainsi $ab$ \`a $a\otimes b$. Soient 
alors $(a_{i\in I})$ (resp. $(b_{i\in I})$) une base
de $A$ (resp. $B$), index\'ees par un ensemble $I$, duales pour $\varphi$, 
{\em i.e.} $\varphi(b_j\otimes a_i)=\delta_{ij}$. Alors l'\'el\'ement
$$\sum_ia_i\otimes b_i$$
est une $R$-matrice universelle pour ${\mathcal D} (A,B,\varphi)$, munissant 
cette 
superalg\`ebre de Hopf d'une structure de superalg\`ebre de Hopf tress\'ee.

\subsection{La superalg\`ebre de Lie $D_x$}

Fixons un nombre complexe $x\not =0,-1$. Rappelons que la superalg\`ebre de Lie
$D_x$ introduite par Kac ({\em cf.} \cite{kac}, o\`u elle est not\'ee
$D(2,1,x)$) admet comme matrice de Cartan la matrice
$$
A=(a_{ij})_{\scriptscriptstyle 1\leq i,j\leq 3}=\left(
\begin{array}{ccc}
0 & 1 & x \\
-1 & 2 & 0 \\
-1 & 0 & 2
\end{array}\right),
$$
et qu'elle a trois racines simples $\alpha_1$ (impaire), $\alpha_2$ et
$\alpha_3$ (paires). Elle admet \'egalement quatre racines positives non simples
qui sont $\alpha_1+\alpha_2$, $\alpha_1+\alpha_3$, $\alpha_1+\alpha_2+\alpha_3$
(impaires) et $2\alpha_1+\alpha_2+\alpha_3$ (paire). Les espaces de racine 
correspondant sont tous de dimension 1. Nous posons
\begin{equation}\label{racines}
\begin{split}
\beta_1=\alpha_3,\quad   \beta_2=\alpha_1+\alpha_3,&\quad  \beta_3=
\alpha_1+\alpha_2+\alpha_3,\\
\beta_4=2\alpha_1+\alpha_2+\alpha_3,\quad  \beta_5=\alpha_1&,\quad
    \beta_6=\alpha_1+\alpha_2,\quad    \beta_7=\alpha_2,
\end{split}\end{equation}
et nous notons $Q=\Z\alpha_1\oplus\Z\alpha_2\oplus\Z\alpha_3$ le {\em r\'eseau
des racines}. Posons
$$
D=(d_i)_{\scriptscriptstyle 1\leq i\leq 3}=\left(
\begin{array}{ccc}
-1 & 0 & 0 \\
0 & 1 & 0 \\
0 & 0 & x 
\end{array}\right),$$ 
de telle fa\c{c}on que la matrice produit
$$
{\overline{A}}=({\overline{a}}_{ij})_{\scriptscriptstyle 1\leq i,j\leq 3}=DA=
\left(\begin{array}{ccc}
0 & -1 & -x \\
-1 & 2 & 0 \\
-x & 0 & 2x
\end{array}\right)
$$
soit sym\'etrique.

Nous aurons besoin du lemme suivant dont la d\'emonstration est laiss\'ee au
lecteur.

\begin{lemm}\label{lem5}
Soient $m,s$ des entiers tels que $1\leq s<m\leq 7$. Alors $n_1\beta_1+\cdots 
+n_s\beta_s\not =n_m\beta_m$ pour tous les \'el\'ements $(n_1,\dots ,n_s,n_m)$ 
non nuls de $\N^{s+1}$, o\`u les $\beta_i$ sont d\'efinis par \eqref{racines}.
\end{lemm}

\subsection{La quantification $\hdx$}\label{hdx}

Nous rappelons ici la d\'efinition de la quantification $\hdx$ de $D_x$, 
puis nous 
d\'efinirons les vecteurs de racine de $\hdx$, {\em cf.} \cite{Zou}. On pourra
\'egalement consulter \cite{CK,FLV,Y96,Z92,ZG91} pour des g\'en\'eralit\'es 
sur les d\'efinitions des supergroupes quantiques.

Nous noterons $\C [[h]]$ l'anneau des s\'eries formelles \`a une 
ind\'etermin\'ee. 
C'est un superespace dont le partie impaire est r\'eduite \`a $0$. Un 
$\C [[h]]$-module $M$ est muni de la topologie $h$-adique pour laquelle une base
de voisinages de $0$ est donn\'ee par la famille $(h^nM)_{n\geq 0}$. Le produit
tensoriel topologique $M\hat{\otimes} N$ de deux $\C [[h]]$-modules est d\'efini 
par
$$M\hat{\otimes} N=\lim_{\leftarrow}\big( M/h^n\! M\otimes_\C N/h^n\! N\big) .$$
Nous renvoyons le lecteur \`a \cite{K} pour
les d\'etails sur la topologie $h$-adique dans le cas g\'en\'eral, et \`a 
\cite{Y94a} dans le cas supergradu\'e. Suivant~\cite{S93}, nous d\'efinissons
 $\hdx$ comme la superalg\`ebre topologiquement engendr\'ee 
par $E_i$, $F_i$, $H_i$, $i=1,2,3$ et les relations (pour tous $i,j=1,2,3$)
\begin{gather}
H_iH_j=H_jH_i, \label{rel1}\\
[H_i, E_j]=a_{ij}E_j, \label{rel2}\\
[H_i,F_j]=-a_{ij}F_j, \label{rel3}\tag{\ref{rel2}$'$}\\
[E_i,F_j]=\delta_{ij}\frac{K_i-K_i^{-1}}{q_i-q_i^{-1}},\notag\\
E_1^2=F_1^2=0,\quad [E_2,E_3]=[F_2,F_3]=0,\notag\\
E_i^2E_1-(q_i+q_i^{-1})E_iE_1E_i+E_1E_i^2=0\text{ si } i=2,3,\notag\\
F_i^2F_1-(q_i+q_i^{-1})F_iF_1F_i+F_1F_i^2=0\text{ si } i=2,3,\notag
\end{gather}
o\`u $\delta_{ij}$ est le symbole de Kronecker, o\`u on a pos\'e
$$q={e}^{h/2},\quad q_i=q^{d_i},\quad K_i=e^{hd_iH_i/2},$$
o\`u $[a,b]=ab-(-1)^{|a||b|}ba$ d\'esigne le supercommutateur et o\`u tous les
g\'en\'erateurs $H_i,E_i,F_i$ sont pairs, sauf $E_1$ et $F_1$ qui sont impairs.
La superalg\`ebre $\hdx$ est une superalg\`ebre de Hopf topologique 
({\em cf.} \cite{Y94a}) dont le coproduit $\Delta$, la co\"unit\'e $\varepsilon$
et l'antipode $S$  sont d\'efinies  par 
\begin{gather}
\Delta (H_i)=H_i\otimes 1+1\otimes H_i,\quad \varepsilon (H_i)=0,\quad 
S(H_i)=-H_i, \label{relbis1}\\
\Delta (E_i)=E_i\otimes 1+K_i\otimes E_i,\quad \varepsilon (E_i)=0,\quad
S(E_i)=-K_i^{-1}E_i, \label{relbis2}\\
\Delta (F_i)=F_i\otimes K_i^{-1}+1\otimes F_i, \quad \varepsilon (F_i)=0,\quad
S(F_i)=-F_iK_i, \label{relbis3}\tag{\ref{relbis2}$'$}
\end{gather}
pour tout $i=1,2,3$. Nous avons \'egalement
$$\Delta (K_i)=K_i\otimes K_i,\quad S(K_i)=K_i^{-1},\quad \varepsilon (K_i)=1.$$

Nous d\'efinissons maintenant l'analogue des vecteurs de racine dans $\hdx$ en 
utilisant les
actions adjointes $\ad$ et $\text{ad}_-$ d\'efinies pour des 
\'el\'ements 
$t,y$ homog\`enes par 
$$\begin{array}{c}
\text{\rm ad}_+t(y)={\displaystyle \sum_{(t)}}(-1)^{|y||t_{(2)}|}\, t_{(1)}y
S(t_{(2)}),\\
\text{\rm ad}_-t(y)={\displaystyle \sum_{(t)}}(-1)^{|y||t_{(1)}|+|t_{(1)}|
|t_{(2)}|}\, t_{(2)}
yS^{-1}(t_{(1)}),
\end{array}$$
(on pourra par exemple consulter \cite{Zou}).
Posons $E_{\beta_1}=E_3$,   $E_{\beta_5}=E_1$,   $E_{\beta_7}=E_2$,
$F_{\beta_1}=F_3$,  $F_{\beta_5}=F_1$,  $F_{\beta_7}=F_2$
et 
\begin{equation}\label{vecteursracines}\begin{array}{ccc}
E_{\beta_2}=\text{\rm ad}_+E_{\beta_5}(E_{\beta_1}), &
E_{\beta_6}=\text{\rm ad}_+E_{\beta_7}(E_{\beta_5}), &
E_{\beta_3}=\text{\rm ad}_+E_{\beta_7}(E_{\beta_2}), \\
E_{\beta_4}=\text{\rm ad}_+E_{\beta_5}(E_{\beta_3}), & &
F_{\beta_2}=\text{\rm ad}_-F_{\beta_5}(F_{\beta_1}),  \\
F_{\beta_6}=\text{\rm ad}_-F_{\beta_7}(F_{\beta_5}),&
F_{\beta_3}=\text{\rm ad}_-F_{\beta_7}(F_{\beta_2}), &
F_{\beta_4}=\text{\rm ad}_-F_{\beta_5}(F_{\beta_3}). \\
\end{array}\end{equation}

\subsection{R\'esultats principaux}\label{resultats}    

Nous \'enon\c{c}ons maintenant nos r\'esultats principaux. Pour des r\'esultats
analogues pour d'autres supergroupes quantiques, on consultera 
\cite{CK,Y94a,Z92,ZG91}. Nous avons besoin des notations suivantes. Posons
\begin{equation}\label{ei_et_phii}
\begin{gathered}
c_2=c_3=c_5=c_6=0,\quad
c_1=2x,\quad c_4=-2-2x,\quad c_7=2,\\
\varphi_1=-\frac{1}{q^{x}-q^{-x}},\quad \varphi_2=\frac{q^{-x}}{q-q^{-1}},\quad 
\varphi_3=-\frac{q^{-1-x}}{q-q^{-1}},\quad
\varphi_6=\frac{q^{-1}}{q-q^{-1}}, \\ 
\varphi_4=-q^{-2-2x}\frac{q^{1+x}-q^{-1-x}}{(q-q^{-1})^2},\quad 
\varphi_5=-\frac{1}{q-q^{-1}}, \quad \varphi_7=-\frac{1}{q-q^{-1}},\\ 
(n)_i=\begin{cases}\frac{q^{nc_i}-1}{q^{c_i}-1} \text{si $c_i\not =0$,}\\ 
1\text{sinon,}\end{cases},\quad
(n)_i!=(n)_i(n-1)_i\cdots (1)_i,\quad n\in \N ,\,1\leq i\leq 7, \\ 
\exp_i(a\otimes b)=\sum_{n\geq 0}\frac{a^n\otimes b^n}{(n)_i!},\quad a,b\in
\hdx ,\quad i=1,\dots ,7,\\
\exp(a\otimes b)=\sum_{n\geq 0}\frac{a^n\otimes b^n}{n!},\quad a,b\in
\hdx ,
\end{gathered}\end{equation}
o\`u $q^x=e^{xh/2}$. Pour $x\not =0,-1$, on pose 
\begin{equation}\label{CetB}
B=(b_{ij})_{\scriptscriptstyle 1\leq i,j\leq 3}=-\frac{1}{2(1+x)}
\left(\begin{array}{ccc}
-4&2&2x\\
2&x&-x\\
2x&-x&x
\end{array}\right).\end{equation}
La matrice $B$ est l'inverse de la matrice $(-a_{ij}/d_j)_
{\scriptscriptstyle 1\leq i,j\leq 3}$.

\begin{theo}\label{theoR}
L'\'el\'ement
\begin{multline*}
R=\exp_1\left(\frac{E_{\beta_1}\otimes F_{\beta_1}}{\varphi_1}\right)\cdots 
\exp_7\left(\frac{E_{\beta_7}\otimes F_{\beta_7}}{\varphi_7}\right)
\exp\left(\frac{h}{2}\Big(\sum_{i,j}b_{ij}H_i\otimes H_j\Big)\right)\\
\in\hdx\hat{\otimes}\hdx
\end{multline*}
est une $R$-matrice universelle pour $\hdx$.
\end{theo}

Nous d\'emontrerons ce th\'eor\`eme au \S\ref{demo_dx}. Posons $R=\sum s_i\otimes 
t_i$ et $u=\sum (-1)^{|s_i|}S(t_i)s_i$. Comme $R\equiv 1\otimes 1\mod h$,
alors $uS(u)\equiv 1\mod h$ admet une
unique racine carr\'ee $\theta\equiv 1\mod h$. Il est bien connu que $\theta$
munit $\hdx$ d'une structure de superalg\`ebre ruban\'ee, {\em cf.} \cite{K,Z95}.
En cons\'equence, la cat\'egorie des $\hdx$-modules topologiques est une cat\'egorie 
ruban\'ee au sens de Turaev (on pourra consulter \cite{Tu2,KRT,K} pour la 
d\'efinition). Pour tout couple
$(V,W)$ de $\hdx$-modules et pour tout $v\in V$, $w\in W$, le {\em tressage} 
$c_{\sst V,W}:V\otimes W\rightarrow W\otimes V$ et le {\em twist} $\theta_{\sst 
V}:V\rightarrow V$ sont donn\'es par
les formules
\begin{equation}\label{twist_et_tressage}
\begin{gathered}
c_{\sst V,W}(v\otimes w)=\tau (R(v\otimes w)),\\
\theta_{\sst V}(v)=\theta^{-1}v.
\end{gathered}
\end{equation}
Au~\S\ref{par:module}, nous introduirons un $\hdx$-supermodule $M$ de rang six.
Le coloriage d'un entrelacs parall\'elis\'e et orient\'e $L$ dans ${\mathbb R}^3$ par 
$M$ fournit un invariant d'isotopie ${\mathcal I}_{L}$ de tels entrelacs.  Avant
d'\'enoncer le th\'eor\`eme~\ref{theo2}, nous rappelons
la d\'efinition de polyn\^ome de Dubrovnik ({\em cf. } \cite{kauffman}). 
C'est un invariant $\Lambda(a,z)$
d'isotopie d'entrelacs parall\'elis\'es non orient\'es \`a valeurs dans
$\Z[a^{\pm 1},z^{\pm 1}]$ qui v\'erifie les relations d'\'echeveaux suivantes :

\begin{gather*}
\Lambda_{\dessin{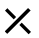}}(a,z)-\Lambda_{\dessin{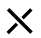}}(a,z)\, =\,
z\big(\Lambda_{\dessin{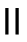}}(a,z)-\Lambda_{\dessin{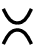}}(a,z)\big),\\
\Lambda_{\dessin{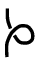}}(a,z)\, =\,a\Lambda_{\dessin{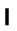}}(a,z),\quad
\Lambda_{\dessin{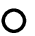}}(a,z)=1.
\end{gather*}

Ces \'egalit\'es font intervenir des entrelacs qui sont \'egaux
sauf dans un disque de ${\mathbb R}^2$ o\`u ils sont comme sur les dessins. 

\begin{theo}\label{theo2}
Pour tout entrelacs parall\'elis\'e orient\'e $L$ dans $\R^3$, l'invariant ${\mathcal I}_L$ 
induit par le supermodule $M$ est li\'e au polyn\^ome de Dubrovnik par
$${\mathcal I}_L=2\Lambda_L(-q^{-1},q-q^{-1}).$$
\end{theo} 
Nous d\'emontrerons ce th\'eor\`eme au \S\ref{demo2}.


\section{Construction d'un double quantique g\'en\'eralis\'e}
 
Dans ce paragraphe, nous construisons une superalg\`ebre de Hopf ${\mathcal D}$ en 
utilisant la th\'eorie du double quantique de Drinfeld ({\em cf.} par exemple 
\cite{KRT}, ainsi que le \S\ref{rappelshopf}). 
Nous utilisons les r\'esultats de \cite{ZGB,Y94a} qui d\'emontrent la validit\'e de 
cette th\'eorie dans le cas supergradu\'e. On pourra \'egalement consulter
\cite{tanisaki}.

\subsection{Pr\'eliminaires techniques}\label{preliminaire}
Dans ce paragraphe, nous d\'efinissons un accouplement de Hopf et nous 
d\'emontrons des lemmes techniques dont nous aurons besoin par la suite.
\begin{enumerate}
\item Nous d\'efinissons $\tilde{U}_+$ comme la superalg\`ebre de Hopf 
topologiquement engendr\'ee sur 
$\C [[h]]$ par
$E_i$ et $H_i$, $(i=1,2,3)$ et les relations~\eqref{rel1} et~\eqref{rel2}.
Le coproduit, la co\"unit\'e et l'antipode sont d\'efinies par les relations
\eqref{relbis1} et~\eqref{relbis2}. Nous d\'efinissons \'egalement les \'el\'ements
$E_{\beta_i}$, $i=1,\dots ,7$ par les relations~\eqref{vecteursracines}.
\item Nous d\'efinissons $\tilde{U}_-$ comme la superalg\`ebre de Hopf 
topologiquement engendr\'ee par
$F_i$ et $H'_i$, $(i=1,2,3)$ et les relations~\eqref{rel1} et~\eqref{rel3}, o\`u
on remplace $H_i$ par $H'_i$.
Le coproduit, la co\"unit\'e et l'antipode sont d\'efinies par les relations
\eqref{relbis1} et~\eqref{relbis3}. Nous d\'efinissons \'egalement les \'el\'ements
$F_{\beta_i}$, $i=1,\dots ,7$ par les relations~\eqref{vecteursracines}. 
\end{enumerate}
Dans tout la suite, nous noterons $q={e}^{h/2}$ et $q_i={e}^{hd_i/2}$.

\noindent {\bf Remarque.} Nous utilisons les m\^emes notations pour les 
g\'en\'erateurs de $\hdx$  que pour ceux de $\tilde{U}_+$, $\tilde{U}_-$. Cet abus 
est justifi\'e {\em a posteriori} par la proposition \ref{isomorphes}.

Notons $\C ((h))$ le corps des fractions de $\C [[h]]$. Posons
\begin{gather}
\label{phi1}
\varphi (F_j\otimes E_i)=-(-1)^{|E_i||F_j|}\frac{\delta_{ij}}{q_i-q_i^{-1}},\\
\label{phi2}
\varphi (H'_j\otimes E_i)=\varphi (F_j\otimes H_i)=0,\\
\label{phi3}
\varphi (H'_j\otimes H_i)=-\frac{2a_{ij}}{hd_j},
\end{gather}
pour $1\leq i,j\leq 3$. D'apr\`es le lemme 3.4 de \cite{KRT} (ou 
\cite{tanisaki}, prop. 2.1.1) qui se g\'en\'eralise 
au cas gradu\'e, les formules \eqref{phi1}-\eqref{phi3} d\'efinissent un 
accouplement de Hopf $\varphi :\tilde{U}_-\otimes\tilde{U}_+\rightarrow
\C ((h))$.

Par souci de concision, nous ne donnons pas les d\'emonstrations des 
lemmes \ref{lem2.1} \`a \ref{iso}, {\em cf.} \cite{moi}.

\begin{lemm}\label{lem2.1}
L'accouplement $\varphi$ v\'erifie,
pour tous $i,j=1,2,3$,
\begin{gather*}
\varphi (K'_j\otimes E_i)=\varphi (K^{'-1}_j\otimes E_i)=
\varphi (F_j\otimes K_i)=\varphi (F_j\otimes K^{-1}_i)=0,\\
\varphi (K'_j\otimes K_i)=\varphi (K^{'-1}_j\otimes K^{-1}_i)=
q^{-\overline{a}_{ij}},\\
\varphi (K'_j\otimes K^{-1}_i)=\varphi (K^{'-1}_j\otimes K_i)=
q^{\overline{a}_{ij}}.
\end{gather*}
\end{lemm}

Soit $\tilde{V}_+$ (resp. $\tilde{V}_-$) la sous-superalg\`ebre de 
$\tilde{U}_+$ (resp. de $\tilde{U}_-$)
engendr\'ee par les $E_i$ (resp. $F_i$), $i=1,2,3$.
La superalg\`ebre $\tilde{V}_+$ (resp. $\tilde{V}_-$) admet une $Q$-graduation 
avec $E_i$ de degr\'e $\alpha_i$ (resp. $F_i$ de degr\'e $-\alpha_i$).
On d\'efinit \'egalement $\tilde{U}_0$ (resp. $\tilde{U}'_0$) comme la 
sous-superalg\`ebre de Hopf de $\tilde{U}_+$ (resp. de $\tilde{U}_-$) 
topologiquement engendr\'ee par les $H_i$ (resp. $H'_i$), $i=1,2,3$.

\begin{lemm}\label{lem2_dx}
Pour tout $H\in \tilde{U}_0$, $F\in \tilde{V}_-$, $F\not =1$, on a
$\varphi (F\otimes H)=0$.
De m\^eme, pour tout $H'\in \tilde{U}'_0$, $E\in \tilde{V}_+$, $E\not =1$, on a
$\varphi (H'\otimes E)=0$.
\end{lemm}

\begin{lemm}\label{lem1}
Si $E\in \tilde{V}_+$, $F\in \tilde{V}_-$, $H\in \tilde{U}_0$, $H'\in 
\tilde{U}'_0$, alors
$$\varphi (FH'\otimes EH)=\varphi (F\otimes E)\varphi (H'\otimes H).$$
\end{lemm}

\begin{lemm}\label{lem3}
Si $i,j_1,\dots ,j_r\in\{1,2,3\}$, et si $r\not =1$ ou $j_1\not =i$, alors
$$\varphi (F_{j_1}\cdots F_{j_r}\otimes E_i)=
\varphi (F_i\otimes E_{j_1}\cdots E_{j_r})=0.$$
\end{lemm}

\begin{lemm}\label{lem4_dx}
Soient $E\in \tilde{V}_+$ de degr\'e $\alpha\in Q$ et
$F\in \tilde{V}_-$ de degr\'e $\beta\in Q$. Si $\alpha +
\beta \not
=0$, alors $\varphi (F\otimes E)=0$.
\end{lemm}

\begin{lemm}\label{lem6}
Soient $n_1,n_4,n_7\in\N$ et $n_2,n_3,n_5,n_6\in\{0,1\}$. Soit $s$ un entier 
compris entre $1$ et $7$. Alors, si $n_i\not =0$ pour $1\leq i<s$ ou $n_s\not
=1$, alors
$$\varphi (\f{s}\otimes \e{1}^{n_1}\cdots \e{s}^{n_s})=
\varphi (\f{1}^{n_1}\cdots \f{s}^{n_s}\otimes \e{s})=0.$$
\end{lemm}

\begin{prop}\label{iso}
Le morphisme de superalg\`ebres $\psi :\tilde{U}_+\rightarrow\tilde{U}_-^
{\text{cop}}$ 
d\'efini sur les g\'en\'erateurs par
$$\psi (E_i)=F_i,\quad \psi (H_i)=-H_i'$$
est un isomorphisme de superalg\`ebres de Hopf qui v\'erifie
$\psi (E_{\beta_i})=F_{\beta_i}$
pour tout $i=1,\dots ,7$.
\end{prop}

Nous d\'efinissons $\tilde{I}_+\subset\tilde{U}_+$ et
$\tilde{I}_-\subset\tilde{U}_-$ comme les annulateurs de $\varphi$ :
\begin{gather*}
\tilde{I}_+=\{E\in\tilde{U}_+,\,\varphi (F\otimes E)=0,\;\forall\;
F\in\tilde{U}_-\},\\
\tilde{I}_-=\{F\in\tilde{U}_-,\,\varphi (F\otimes E)=0,\;\forall\;
E\in\tilde{U}_+\}.
\end{gather*}

\begin{prop}\label{annulateur}
Pour $i=2,3$, les \'el\'ements suivants sont dans $\tilde{I}_+$ (resp. dans 
$\tilde{I}_-$) :
\begin{gather*}
E_1^2,\quad E_2E_3-E_3E_2,\quad E_i^2E_1-(q_i+q_i^{-1})E_iE_1E_i+E_1E_i^2
 \\
\big(\text{resp. } F_1^2,\quad F_2F_3-F_3F_2,\quad F_i^2F_1-(q_i+q_i^{-1})
F_iF_1F_i+
F_1F_i^2 \big).
\end{gather*}
\end{prop}
\begin{proof} Pour montrer que $X\in \tilde{I}_+$, on \'etablit tout d'abord que
$$\Delta (X)=X\otimes 1+K\otimes X,$$ o\`u $K\in \tilde{U}_0$. Ce calcul permet 
alors d'\'ecrire 
$$\varphi (YZ\otimes X)=(-1)^{\xi_1}\varphi (Y\otimes X)\varphi (Y\otimes 1)
+(-1)^{\xi_2}\varphi (Y\otimes K)\varphi (Z\otimes X),$$
qui est nul si $\varphi (Y\otimes X)=\varphi (Z\otimes X)=0$. On montre alors
que $\varphi (Y\otimes X)=0$ pour tout g\'en\'erateur $Y\in\tilde{U}_-$, 
{\em i.e.}
pour $Y=H'_j$ et $Y=F_j$. Ce dernier point d\'ecoule des lemmes \ref{lem2_dx} et 
\ref{lem3}. Nous calculons donc uniquement les coproduits des \'el\'ements 
consid\'er\'es. Pour cela, nous rappelons que nous utilisons la formule
\eqref{produit} pour le produit de deux \'el\'ements de $\tilde{U}_+
\hat{\otimes}
\tilde{U}_+$. Nous traitons uniquement le cas de
$z=E_i^2E_1-(q_i+q_i^{-1})E_iE_1E_i+E_1E_i^2$, $i=2$ ou $i=3$.

Soient \mbox{$1\leq n,m,p\leq 3$} des entiers tels qu'un seul d'entre eux
soit \'egal \`a $1$. On a alors les deux \'egalit\'es suivantes.
\begin{multline*}
\Delta (E_nE_m)=E_nE_m\otimes 1+E_nK_m\otimes E_m+K_nE_m\otimes E_n+
K_nK_m\otimes E_nE_m\\
=E_nE_m\otimes 1+E_nK_m\otimes E_m+q^{{\overline{a}}_{nm}}E_mK_n\otimes E_n+
K_nK_m\otimes E_nE_m,\vspace*{-5cm}
\end{multline*}
\vspace*{-.7cm}
\begin{multline*}
\Delta (E_nE_mE_p)=E_nE_mE_p\otimes 1+E_nK_mE_p\otimes E_m+
q^{{\overline{a}}_{nm}}
E_mK_nE_p\otimes E_n \\
+K_nK_mE_p\otimes E_nE_m
+E_nE_mK_p\otimes E_p+E_nK_mK_p\otimes E_mE_p \\
+q^{{\overline{a}}_{nm}}E_mK_nK_p\otimes E_nE_p+K_nK_mK_p\otimes E_nE_mE_p \\
=E_nE_mE_p\otimes 1+q^{{\overline{a}}_{mp}}E_nE_pK_m\otimes E_m+
q^{{\overline{a}}_{nm}+{\overline{a}}_{np}}
E_mE_pK_n\otimes 
E_n\\
+q^{{\overline{a}}_{mp}+{\overline{a}}_{np}}E_pK_nK_m\otimes E_nE_m
\end{multline*}
On en d\'eduit
\begin{multline*}
\Delta (E_i^2E_1)=E_i^2E_1\otimes 1+K_i^2K_1\otimes E_i^2E_1
+q^{{\overline{a}}_{i1}}(1+q^{{\overline{a}}_{ii}})E_iE_1K_i\otimes E_i \\
+E_i^2K_1\otimes E_1
+q^{2{\overline{a}}_{i1}}E_1K_i^2\otimes E_i^2+(1+q^{{\overline{a}}_{ii}})
E_iK_iK_1\otimes E_iE_1,
\end{multline*}
\vspace*{-.7cm}
\begin{multline*}
\Delta (E_1E_i^2)=E_1E_i^2\otimes 1+K_1K_i^2\otimes E_1E_i^2+
E_1K_i^2\otimes E_i^2 \\
+(1+q^{{\overline{a}}_{ii}})E_1E_iK_i\otimes E_i+
q^{2{\overline{a}}_{1i}}E_i^2K_1\otimes E_1+q^{{\overline{a}}_{1i}}(1+
q^{{\overline{a}}_{ii}})E_iK_1K_i\otimes 
E_1E_i,
\end{multline*}
\vspace*{-.7cm}
\begin{multline*}
\Delta (E_iE_1E_i)=E_iE_1E_i\otimes 1+E_iE_1K_i\otimes E_i\\
+q^{{\overline{a}}_{1i}}E_i^2K_1\otimes E_1+E_iK_1K_i\otimes E_1E_i
+q^{{\overline{a}}{i1}+{\overline{a}}_{ii}}E_1E_iK_i\otimes E_i+
q^{{\overline{a}}_{i1}}E_1K_i^2\otimes E_i^2\\+
q^{{\overline{a}}{1i}+{\overline{a}}_{ii}}E_iK_iK_1\otimes E_iE_1+
K_i^2K_1\otimes E_iE_1E_i.
\end{multline*}
Ceci prouve que $\Delta (z)=z\otimes 1+K_1K_i^2\otimes z+X$, 
o\`u $X$ est un 
\'el\'ement dont il reste \`a montrer qu'il est nul. Or,
\begin{multline*}
X=(1+q^{2{\overline{a}}_{1i}}-q^{{\overline{a}}_{1i}}(q_i+q_i^{-1}))E_i^2K_1\otimes E_1
+(1+q^{2{\overline{a}}_{i1}}-(q_i+q_i^{-1})q^{{\overline{a}}_{i1}})E_1K_i^2\otimes E_i^2\\
+(q^{{\overline{a}}_{i1}}(1+q^{{\overline{a}}_{ii}})-(q_i+q_i^{-1}))E_iE_1K_i\otimes E_i \\
+(1+q^{{\overline{a}}_{ii}}-q^{{\overline{a}}_{i1}+{\overline{a}}_{ii}}(q_i+q_i^{-1}))E_1E_iK_i\otimes E_i\\
+(1+q^{{\overline{a}}_{ii}}-q^{{\overline{a}}_{1i}+{\overline{a}}_{ii}}(q_i+q_i^{-1}))E_iK_iK_1\otimes E_iE_1\\
+(q^{{\overline{a}}_{1i}}(1+q^{{\overline{a}}_{ii}})-(q_i+q_i^{-1}))E_iK_1K_i\otimes E_1E_i,
\end{multline*}
et pour $i=2,3$, on a ${\overline{a}}_{i1}={\overline{a}}_{1i}=d_ia_{i1}=-d_i$ et ${\overline{a}}_{ii}=2d_i$.
On en d\'eduit
\begin{gather*}
1+q^{2{\overline{a}}_{1i}}-q^{{\overline{a}}_{1i}}(q_i+q_i^{-1})=1+q^{-2d_i}-q^{-d_i}
(q_i+q_i^{-1})=0,\\
q^{{\overline{a}}_{i1}}(1+q^{{\overline{a}}_{ii}})-(q_i+q_i^{-1})=q^{-d_i}(1+q^{2d_i})-(q^{d_i}
+q^{-d_i})=0,\\
1+q^{{\overline{a}}_{ii}}-q^{{\overline{a}}_{i1}+{\overline{a}}_{ii}}(q_i+q_i^{-1})=1+q^{2d_i}-q^{-d_i}(q^{d_i}
+q^{-d_i})=0,
\end{gather*}
d'o\`u $X=0$ et la proposition est d\'emontr\'ee en ce qui concerne $\tilde{I}_+$.
Pour montrer que $Y\in\tilde{I}_-$, il suffit d'\'etablir que 
$\Delta (Y)=1\otimes Y+Y\otimes K'$, o\`u $K'\in\tilde{U}'_0$, ce qui d\'ecoule
des calculs pr\'ec\'edents en appliquant l'isomorphisme $\psi$.\end{proof}

On note $I_+$ (resp. $I_-$) le sous-id\'eal de Hopf de $\tilde{I}_+$ (resp.
$\tilde{I}_-$) engendr\'e par
$$\begin{array}{rl}
& E_1^2,\, E_2E_3-E_3E_2,\, E_i^2E_1-(q_i+q_i^{-1})E_iE_1E_i+E_1E_i^2\, 
(i=2,3) \\
\text{(resp. par}&F_1^2,\, F_2F_3-F_3F_2,\, F_i^2F_1-(q_i+q_i^{-1})F_iF_1F_i+
F_1F_i^2\, (i=2,3)).
\end{array}$$
On pose alors
$$U_+=\tilde{U}_+/I_+,\quad U_-=\tilde{U}_-/I_-.$$

\begin{coro}
L'accouplement $\varphi$ induit un accouplement sur
$U_-\hat{\otimes} U_+$ que nous noterons encore $\varphi$. De plus, l'isomorphisme
$\psi$ de la proposition \ref{iso} induit un isomorphisme de superalg\`ebres
de Hopf entre $U_+$ et 
$U_-^{\text{cop}}$, que nous noterons encore $\psi$.
\end{coro}

\subsection{La superalg\`ebre de Hopf ${\mathcal D}$}\label{D}
Dans ce paragraphe, nous construisons une superalg\`ebre de Hopf tress\'ee 
${\mathcal D}$ 
par la m\'ethode du double quantique ({\em cf.} \S\ref{rappelshopf}) et nous
\'etablissons un lien avec $\hdx$.

En paraphrasant les d\'efinitions de $\tilde{V}_+,\tilde{V}_-,\tilde{U}_0,
\tilde{U}'_0$ du \S\ref{preliminaire} on d\'efinit $V_+\subset U_+$, 
$V_-\subset U_-$, $U_0\subset U_+$, $U'_0\subset U_-$.
Les r\'esultats du lemme \ref{lem2_dx} jusqu'\`a la proposition \ref{iso} 
restent
valables pour $\varphi:U_-\otimes U_+\rightarrow\C [[h]]$, en rempla\c{c}ant
$\tilde{V}_+$ (resp. $\tilde{V}_-$, $\tilde{U}_0$, $\tilde{U}'_0$) par $V_+$ 
(resp. $V_-$, $U_0$, $U'_0$). Nous notons ${\mathcal D} =
{\mathcal D}(U_+,U_-,\varphi)$ le double
quantique de $U_+$ et $U_-$ construit \`a partir de l'accouplement $\varphi$
({\em cf.} le \S\ref{rappelshopf} pour la construction).

Nous avons besoin de la proposition suivante pour \'etablir un lien entre 
${\mathcal D}$ et $\hdx$. Nous laissons la d\'emonstration au lecteur.

\begin{prop}\label{relations}
Les relations suivantes sont v\'erifi\'ees dans ${\mathcal D}$:
\begin{gather*}
[H'_i, E_j]=a_{ij}E_j,\quad [H_i, F_j]=-a_{ij}F_j, \\
[E_i,F_j]=\delta_{ij}\frac{K_i-K_i^{'-1}}{q_i-q_i^{-1}}.
\end{gather*}
\end{prop}

D\'efinissons la superalg\`ebre ${\overline{\mathcal D}}$ comme le quotient de 
${\mathcal D}$ par l'id\'eal de Hopf engendr\'e par $H_i-H_i'$, pour $i=1,2,3$. 
Posons
$\pi :{\mathcal D}\rightarrow{\overline{\mathcal D}}$ la projection canonique. 
D\'efinissons 
$\iota :\hdx\rightarrow{\mathcal D}$ par
$$\iota (E_j)=E_j,\quad \iota (F_j)=F_j,\quad \iota (H_j)=H_j,$$
pour $j=1,2,3$.

\begin{prop}\label{isomorphes}
Le morphisme $\iota$ est un morphisme de superalg\`ebres de 
Hopf topologiques injectif et la compos\'ee $p\circ\iota$ r\'ealise un 
isomorphisme entre $\hdx$ et ${\overline{\mathcal D}}$.
\end{prop}


\section{$R$-matrice universelle}
Dans ce paragraphe, nous d\'emontrons le th\'eor\`eme \ref{theoR}
en construisant des bases de $U_+$ et $U_-$ duales pour ~$\varphi$. Rappelons
que $D_x$ admet trois racines simples $\alpha_1$, $\alpha_2$ et 
$\alpha_3$ ({\em cf.} \S\ref{hdx}), et que les vecteurs de racine sont d\'efinis
par \eqref{vecteursracines}. Nous notons \'egalement $\e{i}$ et $\f{i}$ les 
images de $\e{i}$ et $\f{i}$ par l'injection $\iota$ de la proposition 
\ref{isomorphes}. 

Pour tout couple $(Y,Z)$ d'\'el\'ements homog\`enes de ${\mathcal D}$ et tout 
$\alpha\in \C[[h]]$, on pose
$$[Y,Z]_{\alpha}=YZ-(-1)^{|Y||Z|}\alpha ZY,$$
et on \'etend cette d\'efinition \`a tout ${\mathcal D}\times {\mathcal D}$ par 
bilin\'earit\'e.

\subsection{Relations dans ${\mathcal D}$}\label{relations_D}
Nous commen\c{c}ons par \'enoncer des relations v\'erifi\'ees dans 
${\mathcal D}$. Les
d\'emonstrations des lemmes~\ref{lem9bis} \`a \ref{ad} sont laiss\'ees au lecteur.
\begin{lemm}\label{lem9bis}
Pour tout $i=1,\dots ,7$, on a
$K_{\beta_i}E_{\beta_i}=q^{c_i}E_{\beta_i}K_{\beta_i}$
\end{lemm}

\begin{lemm}\label{valeur_racines} On a
\begin{gather*}
E_{\beta_2}=[E_1,E_3]_{q^{-x}},\qquad E_{\beta_6}=[E_2,E_1]_{q^{-1}}, \\
E_{\beta_3}=[E_2,E_{\beta_2}]_{q^{-1}}, \qquad
E_{\beta_4}=[E_1,E_{\beta_3}]_{q^{-1-x}}.
\end{gather*}
\end{lemm}

\begin{lemm}\label{ad}
L'action adjointe $\ad:U_+\rightarrow U_+$ v\'erifie pour $i=2$ et $3$
$$\ad^2 E_i(E_1)=0,$$
et pour tout $X\in\, U_+$, on a
$$\ad E_i(XE_j)=\ad E_i(X)\, E_j,\quad \ad E_i(E_jX)=E_j\,\ad E_i(X),$$
si $(i,j)=(2,3)$ ou $(3,2)$.
\end{lemm}

Le lemme suivant a \'et\'e \'etabli par Zou dans \cite{Zou}.
\begin{lemm}\label{carres}
Pour $i=2,3,6$, on a $E_{\beta_i}^2=0$.
\end{lemm}

\begin{prop}\label{commut}
Dans $U_+$, on a les relations de commutations suivantes:
\begin{equation}\label{type1}
\begin{gathered}
 \phantom{a} [\e{7},\e{5}]_{q^{-1}}=\e{6},\quad [\e{7},\e{2}]_{q^{-1}}=\e{3},
\quad [\e{7},\e{1}]=0,\\
[\e{5},\e{3}]_{-q^{-1-x}}=\e{4},\quad [\e{5},\e{1}]_{q^{-x}}=\e{2},
\end{gathered}
\end{equation}
\begin{equation}\label{type2}
\begin{gathered}
\phantom{a} [\e{7},\e{6}]_q=0,\quad [\e{7},\e{3}]_q=0,\quad
[\e{6},\e{5}]_{-q^{-1}}=0,\\ 
 [\e{6},\e{1}]_{q^{-x}}=\e{3},\quad  [\e{5},\e{4}]_{q^{-1-x}}=0,\quad
 [\e{5},\e{2}]_{-q^{-x}}=0,\\
 [\e{4},\e{3}]_{q^{-1-x}}=0,\quad
 [\e{3},\e{2}]_{-q^{-1}}=0,\quad [\e{2},\e{1}]_{q^{x}}=0,\\
\end{gathered}
\end{equation}
\begin{equation}\label{type3}
\begin{gathered}
\phantom{a} [\e{7},\e{4}]=0,\quad [\e{6},\e{4}]_{q^{-1-x}},\quad
[\e{6},\e{3}]_{-q^{-x}}=0,\\
[\e{4},\e{2}]_{q^{-1-x}}=0,\quad
[\e{4},\e{1}]=q^{-1}(q^{1+x}-q^{-1-x})\e{2}\e{3},\\
[\e{3},\e{1}]_{q^{x}}=0,\quad
[\e{6},\e{2}]_{-q^{-1-x}}=-q^{-1-x}(q-q^{-1})\e{3}\e{5}-q^{-1}\e{4}.
\end{gathered}
\end{equation}
\end{prop}
\begin{proof} Les relations~\eqref{type1} sont imm\'ediates, et nous laissons au lecteur
la d\'emonstration des relations~\eqref{type2}. Nous d\'emontrons les 
relations~\eqref{type3}.
Nous utiliserons implicitement les lemmes \ref{valeur_racines},
\ref{ad} et \ref{carres}, et les relations~\eqref{type2}.
\begin{enumerate}
\item On a
\begin{align*}
\e{7}\e{4}&=\e{7}(\e{5}\e{3}+q^{-1-x}\e{3}\e{5}) \\
&=q^{-1}\e{5}\e{7}\e{3}+\e{6}\e{3}+q^{-x}\e{3}\e{7}\e{5}\\
&=\e{5}\e{3}\e{7}+\e{6}\e{3}+q^{-1-x}\e{3}\e{5}\e{7}+q^{-x}\e{3}\e{6}
=\e{4}\e{7},
\end{align*}
o\`u la deuxi\`eme et la troisi\`eme \'egalit\'e proviennent des relations de 
commutations de $\e{7}\e{5}$ et $\e{7}\e{3}$ et la derni\`ere de celle de
$\e{6}\e{3}$.
\item On a
\begin{align*}
\e{6}\e{4}&=\e{6}(\e{5}\e{3}+q^{-1-x}\e{3}\e{5})
=\e{6}\e{5}\e{3}+q^{-1-x}\e{6}\e{3}\e{5}\\
&=-q^{-1}\e{5}\e{6}\e{3}-q^{-1-2x}\e{3}\e{6}\e{5}\\
&=-q^{-1-x}\e{5}\e{3}\e{6}-q^{-2-2x}\e{3}\e{5}\e{6}\\
&=q^{-1-x}\e{4}\e{6},
\end{align*}
o\`u on a utilis\'e les relation de commutation de $\e{6}\e{5}$ et
$\e{6}\e{3}$.
\item On a
$\e{6}\e{3}+q^{-x}\e{3}\e{6}=\e{6}^2E_3-q^{-x}\e{6}E_3\e{6}+q^{-x}\e{6}E_3\e{6}
-q^{-2x}E_3\e{6}^2=0$
car $\e{6}^2=0$ et $\e{3}=\e{6}\e{1}-q^{-x}\e{1}\e{6}$.
\item En rempla\c{c}ant $\e{4}$ par sa valeur, on obtient
\begin{align*}
\e{4}\e{2}&=\e{5}\e{3}\e{2}+q^{-1-x}\e{3}\e{5}\e{2}\\
&=q^{-1-x}\e{2}\e{5}\e{3}+q^{-2-2x}\e{2}\e{3}\e{5}.
\end{align*}
\item On a
\begin{align*}
\e{3}\e{1}-q^{x}\e{1}\e{3}&=\ad E_2(\e{2})E_3-q^{x}E_3\ad E_2(\e{2})\\
&=\ad E_2(\e{2}E_3-q^{x}E_3\e{2})=0
\end{align*}
car $\e{2}\e{1}=q^{x}\e{1}\e{2}$.
\item On remplace \`a nouveau $\e{4}$ par sa valeur donn\'ee au lemme 
\ref{valeur_racines} pour obtenir
$$\e{4}\e{1}=q^x\e{5}\e{1}\e{3}+q^{-1-x}\e{3}\e{5}\e{1}$$
car $\e{3}\e{1}=q^x\e{1}\e{3}$. Comme $\e{5}\e{1}=q^{-x}\e{1}\e{5}+\e{2}$, nous
obtenons
$$\e{4}\e{1}=\e{1}\e{3}\e{5}+q^{x}\e{2}\e{3}+q^{-1-x}\e{1}\e{5}\e{3}
+q^{-1-x}\e{3}\e{2},$$
ce qui donne le r\'esultat voulu avec $\e{3}\e{2}=-q^{-1}\e{2}\e{3}$.
\item En rempla\c{c}ant $\e{6}$ par sa valeur donn\'ee au lemme 
\ref{valeur_racines}, nous obtenons
\begin{align*}
\e{6}\e{2}&=\e{7}\e{5}\e{2}-q^{-1}\e{5}\e{7}\e{2}\\
&=-q^{-x}\e{7}\e{2}\e{5}-q^{-2}\e{5}\e{2}\e{7}-q{-1}\e{5}\e{3}\\
&=-q^{-1-x}\e{2}\e{7}\e{5}-q^{-x}\e{3}\e{5}+q^{-2-x}\e{2}\e{5}\e{7}\\
\phantom{aaaaaaaa}&+q^{-2-x}\e{3}\e{5}-q^{-1}\e{4}\\
&=-q^{-1-x}\e{2}\e{6}-q^{-1-x}(q-q^{-1})\e{3}\e{5}-q^{-1}\e{4},
\end{align*}
o\`u on a utilis\'e les relations de commutation de $\e{5}\e{2}$, $\e{7}\e{2}$
et $\e{5}\e{3}$.
\end{enumerate}
\end{proof}

\begin{prop}\label{coproduits}
Les coproduits des vecteurs de racines non simples sont donn\'es par les formules
suivantes:
\begin{gather*}
\Delta (\e{2})=\e{2}\otimes 1+\k{{2}}\otimes \e{2}+(1-q^{-2x})\e{5}\k{{1}}
\otimes \e{1}, \\
\Delta (\e{6})=\e{6}\otimes 1+\k{{6}}\otimes \e{6}+(1-q^{-2})\e{7}\k{{5}}
\otimes \e{5},
\end{gather*}
\begin{multline*}
\Delta (\e{3})=\e{3}\otimes 1+\k{{3}}\otimes \e{3}+(1-q^{-2x})\e{6}\k{{1}}
\otimes \e{1} \\
 +(1-q^{-2})\e{7}\k{{2}}\otimes \e{2},
\end{multline*}
\begin{multline*}
\Delta (\e{4})=\e{4}\otimes 1+\k{{4}}\otimes \e{4}+(1-q^{-2-2x})\e{5}\k{{3}}
\otimes \e{3} \\
-q^{-1}(1-q^{-2-2x})\e{6}\k{{2}}\otimes \e{2} 
 +(1-q^{-2x})(1-q^{-2-2x})\e{5}\e{6}\k{{1}}\otimes \e{1} \\
 +(1-q^{-2})(1-q^{-2-2x})\e{5}\e{7}\k{{2}}\otimes \e{2}.
\end{multline*}
\end{prop}
\begin{proof} Nous traitons ici uniquement le cas de $\e{4}$ qui est le plus
compliqu\'e. On a
\begin{multline*}
\Delta (E_1)\Delta (\e{3})=(E_1\otimes 1+K_1\otimes E_1)
(\e{3}\otimes 1+\k{{3}}\otimes \e{3} \\ 
+(1-q^{-2x})\e{6}\k{{1}}\otimes \e{1}+(1-q^{-2})\e{7}\k{{2}}\otimes \e{2})\\
=E_1\e{3}\otimes 1+E_1\k{{3}}\otimes \e{3}+(1-q^{-2})E_1\e{7}\k{{2}}\otimes 
\e{2}\\
+(1-q^{-2x})E_1\e{6}\k{{1}}\otimes \e{1}
-q^{-1-x}\e{3}K_1\otimes E_1+K_1\k{{3}}\otimes E_1\e{3}\\
+q^{-1}(1-q^{-2})\e{7}K_1\k{{2}}\otimes E_1\e{2}-q^{-1}(1-q^{-2x})\e{6}K_1
\k{{1}}\otimes E_1\e{1},
\end{multline*}
o\`u nous avons utilis\'e les relations de commutation entre les $K_i$ et les
$E_j$. De m\^eme, on a
\begin{multline*}
\Delta (\e{3})\Delta (E_1)=
\e{3}E_1\otimes 1-q^{-1}E_1\k{{3}}\otimes \e{3}
-q^{-x}(1-q^{-2})\e{7}E_1\k{{2}}
\otimes \e{2}\\
+q^{-x}(1-q^{-2x})\e{6}E_1\k{{1}}
\otimes \e{1}
+\e{3}K_1\otimes E_1+\k{{3}}K_1\otimes \e{3}E_1\\
+(1-q^{-2})\e{7}\k{{2}}K_1\otimes \e{2}E_1
+(1-q^{-2x})\e{6}\k{{1}}K_1\otimes \e{1}E_1,
\end{multline*}
d'o\`u on d\'eduit que
\begin{multline*}
\Delta (\e{4})=
\e{4}\otimes 1+\k{4}\otimes \e{4}+(1-q^{-2x})(E_1\e{6}+q^{-1-2x}
\e{6}E_1)K_3\otimes \e{1}\\
+(1-q^{-2-2x})\e{5}\k{3}\otimes\e{3}
+(1-q^{-2})(E_1E_2-q^{-1-2x}E_2E_1)\k{2}\otimes\e{2}\\
-q^{-1}(1-q^{-2x})\e{6}K_1K_3\otimes (E_1E_3-q^{-x}E_3E_1) \\
+q^{-1}(1-q^{-2-x})E_2K_1\k{2}\otimes (E_1\e{2}+q^{-x}\e{2}E_1).
\end{multline*}
Or, par la proposition \ref{commut}, on a
$E_1\e{6}+q^{-1-2x}\e{6}E_1=(1-q^{-2-2x})\e{5}\e{6}$,
et
$$-q^{-1}(1-q^{-2x})\e{6}K_1K_3\otimes (E_1E_3-q^{-x}E_3E_1)=
-q^{-1}(1-q^{-2x})\e{6}\k{2}\otimes \e{2}.$$
Ensuite, nous calculons $(1-q^{-2})(E_1E_2-q^{-1-2x}E_2E_1)\k{2}\otimes\e{2}$.
En utilisant la proposition \ref{commut}, on obtient
\begin{multline*}
(1-q^{-2})(E_1E_2-q^{-1-2x}E_2E_1)\k{2}\otimes\e{2} \\
=(1-q^{-2})(\e{5}\e{7}-q^{-1-2x}
\e{7}\e{5})\k{2}\otimes\e{2}\\
=(1-q^{-2})(\e{5}\e{7}-q^{-2-2x}\e{5}\e{7}-q^{-1-2x}\e{6})\k{2}\otimes\e{2}.
\end{multline*}
Pour terminer, il reste le terme
$q^{-1}(1-q^{-2-x})E_2K_1\k{2}\otimes (E_1\e{2}+q^{-x}\e{2}E_1).$
Or, d'apr\`es la proposition \ref{commut}, on a
$$E_1\e{2}+q^{-x}\e{2}E_1=\e{5}\e{2}+q^{-x}\e{2}\e{5}=0,$$
ce qui ach\`eve la d\'emonstration de la proposition.\end{proof}

\subsection{D\'emonstration du th\'eor\`eme \ref{theoR}}\label{demo_dx}
Nous donnons maintenant une d\'emonstration du th\'eor\`eme \ref{theoR}.

\begin{lemm}\label{lem9}
Pour tout $i=1,\dots ,7$, on a
$\varphi (F_{\beta_i}\otimes E_{\beta_i})=\varphi_i$,
o\`u $\varphi_i$, $i=1,\dots ,7$ sont d\'efinis par 
~\eqref{ei_et_phii}.
\end{lemm}
\begin{proof} Elle est laiss\'ee au lecteur.\end{proof}

\begin{prop}\label{delta(X)}
Soient $n_1,n_4,n_7\in\N$, $n_2,n_3,n_5,n_6\in\{0,1\}$ et $X=E_{\beta_1}^{n_1}
\cdots E_{\beta_7}^{n_7}\in V_+$. Alors
$$\Delta (X)=(n_s)_s\, \e{1}^{n_1}\cdots \e{s-1}^{n_{s-1}}
\e{s}^{n_s-1}\k{s}\otimes \e{s}+ y\otimes 1+
\sum z\otimes \left(\prod_{i\leq s}\e{i}\right),$$
o\`u $y,z\in U_+$,  $1\leq s\leq 7$ est le plus grand des entiers $i$ tel que 
$n_i\not =0$, et o\`u $\left(\prod_{i\leq s}E_{\beta_i}\right)$ 
d\'esigne un produit de racines $\e{i}$ dans l'ordre croissant des indices, 
avec au moins un des indices $i$ v\'erifiant $i<s$.  
\end{prop}
\begin{proof} D'apr\`es la proposition \ref{coproduits} et les relations 
\eqref{relbis2}, pour $i=1,\dots ,7$, on a
$$\Delta (E_{\beta_i})=E_{\beta_i}\otimes 1+K_{\beta_i}\otimes E_{\beta_i}
+\sum_{\beta_k,\, k<i}\, z\otimes E_{\beta_k}$$
o\`u $z\in U_+$. On en d\'eduit que
$$\Delta (X)=(\e{1}\otimes 1+\k{1}\otimes\e{1})^{n_1}\cdots 
(\e{s}\otimes 1+\k{s}\otimes\e{s}+\sum_{\beta_k,k<s} z\otimes\e{k})^{n_s}.$$
Le terme $y\otimes 1$ de l'\'enonc\'e en d\'ecoule aussit\^ot.
Consid\'erons le terme $z\otimes\e{s}$, $z\in U_+$. 
Pour $k=0,\dots ,n_s-1$, on a (en utilisant le lemme \ref{lem9bis}),
\begin{multline*}
(\e{1}^{n_1}\cdots\e{s-1}^{n_{s-1}}\otimes 1)(\e{s}^k\otimes 1)(\k{s}\otimes
\e{s})(\e{s}^{n_s-k-1}\otimes 1)\\
=\e{1}^{n_1}\cdots\e{s-1}^{n_{s-1}}\e{s}^k\k{s}\e{s}^{n_s-k-1}\otimes\e{s} \\
=q^{c_s(n_s-k-1)}\e{1}^{n_1}\cdots\e{s-1}^{n_{s-1}}\e{s}^{n_s-1}\k{s}\otimes
\e{s},
\end{multline*}
 car aucun signe ne peut appara\^itre dans 
le produit, (si $\e{s}$ est impair,
alors $n_s=1$ d'apr\`es le lemme \ref{carres}). Or,
$$\sum_{k=0}^{n_s-1}q^{kc_s}=\frac{q^{n_sc_s}-1}{q^{c_s}-1}=(n_s)_s,$$
ce qui fournit le premier terme de l'\'enonc\'e.
Il reste \`a calculer le produit
$$\left(\sum_{\beta_k,k<i}z\otimes\e{k}\right)\left(\sum_{\beta_m,m<j}z'\otimes
\e{m}\right),$$
avec $i<j\leq s$. Ce produit donne 
des termes  $z\otimes \e{p}\e{r}$ avec
$p>r$. Or, la proposition \ref{commut} prouve que
$\e{p}\e{r}=\lambda\e{r}\e{p}+\zeta\e{f}\e{g}+\mu\e{t},$
o\`u $\lambda ,\zeta ,\mu\in\C[[h]]$ et $r<f,g,t<p$. En particulier, avec
$p\leq s$, on
obtient le dernier terme de l'\'enonc\'e.\end{proof}

\begin{lemm}\label{phi(F_E)}
Soient $n_i,m_i\in\N$, $i=1,4,7$, $n_j,m_j\in\{0,1\}$, $j=2,3,5,6$. Alors
$$\varphi (F_{\beta_1}^{m_1}\cdots F_{\beta_7}^{m_7}\otimes E_{\beta_1}^{n_1}
\cdots E_{\beta_7}^{n_7})=(-1)^\xi
\prod_{i=1}^7\delta_{n_im_i}\, (n_i)_i!\;\varphi_i^{n_i},$$
o\`u $\xi=n_6(n_2+n_3+n_5)+n_5(n_2+n_3)+n_3n_2$.
\end{lemm}
\begin{proof} Soit $s$ le plus grand des indices $i$ tel que $m_i$ ou $n_i$ soit
non nul. Nous supposerons pour la d\'emonstration que $n_s\geq m_s$. Le cas
$m_s\geq n_s$ se traite de la m\^eme mani\`ere gr\^ace \`a la 
proposition~\ref{iso}.
Posons $X=E_{\beta_1}^{n_1}\cdots E_{\beta_{s-1}}^{n_{s-1}}E_{\beta_s}
^{n_s-1}$ et $Y=F_{\beta_1}^{m_1}\cdots F_{\beta_{s-1}}^{m_{s-1}}F_{\beta_s}^
{m_s-1}$. Si $m_s\not=0$, alors
\begin{multline*}
\varphi (YF_{\beta_s}\otimes X\e{s}) \\
=\sum_{(X\e{s})}(-1)^{|\f{s}||(X\e{s})_{(1)}|}\,
\varphi (Y\otimes (X\e{s})_{(1)})\,\varphi (\f{s}\otimes (X\e{s})_{(2)}).
\end{multline*}
Or, d'apr\`es la proposition \ref{delta(X)}, on a $\varphi (\f{s}\otimes (X\e{s}
)_{(2)})=0$, sauf lorsque $(X\e{s})_{(2)}=\e{s}$. En effet, si $(X\e{s})_
{(2)}=1$, cela d\'ecoule du lemme \ref{lem2_dx}. Sinon,
$(X\e{s})_{(2)}=\prod_{i\leq s}\e{i}$ par la proposition \ref{delta(X)}, et le
lemme \ref{lem6} permet de conclure. On en d\'eduit que 
\begin{align*}
\varphi (YF_{\beta_s}\otimes X\e{s})&=(-1)^{|\f{s}||(X\e{s})_{(1)}|}\,(n_s)_s\,
\varphi_s\,\varphi (Y\otimes X\k{s})\\
&=(-1)^{|\f{s}||(X\e{s})_{(1)}|}\,
(n_s)_s\,\varphi_s\,\varphi (Y\otimes X)
\end{align*}
d'apr\`es le lemme \ref{lem1}.
Or, on a l'\'egalit\'e des degr\'es 
$$|(X\e{s})_{(1)}|=|E_{\beta_1}^{n_1}\cdots 
E_{\beta_{s-1}}^{n_{s-1}}E_{\beta_s}^{n_s-1}K_{\beta_s}|.$$ 
Le seul cas qui nous 
int\'eresse est celui o\`u $|\beta_s|=1$, et dans ce cas $n_s-1=0$. On en 
d\'eduit que
$$
|(X\e{s})_{(1)}|=|E_{\beta_1}^{n_1}\cdots E_{\beta_{s-1}}^{n_{s-1}}|
=\sum_{i<s,\, |\beta_i|=1}n_i$$
puisque, lorsque $\beta_i$ est impaire, nous avons $n_i=0$ ou $1$ (lemme
\ref{carres}). En it\'erant cette op\'eration $n_s-m_s$ fois, on obtient
\begin{multline*}
\varphi (Y\f{s}\otimes X\e{s}) \\
=(-1)^{\zeta_s}(n_s)_s\cdots (n_s-m_s+1)_s\,
\varphi^{n_s-m_s}\,\varphi (\f{1}^{m_1}\cdots \f{m_{s-1}}^{m_{s-1}}
\otimes \e{1}^{n_1}\cdots \e{s}^{n_s-m_s}),
\end{multline*}
o\`u $\zeta_s=\sum_{i<s,\, |\beta_i|=1}n_i$. Si $n_s=m_s$, alors
on peut refaire le raisonnement pr\'ec\'edant tant que $n_i=m_i$, et ce jusqu'\`a ce 
que $s=1$. De plus, il appara\^it un signe $(-1)^\xi$ o\`u 
$$\xi=\sum_{i=1}^{7}
\zeta_i=n_6(n_2+n_3+n_5)+n_5(n_2+n_3)+n_3n_2.$$ 
Il reste donc \`a montrer que
$$\varphi (\f{1}^{m_1}\cdots \f{r}^{m_r}\otimes\e{1}^{n_1}\cdots\e{s}^{n_s}
)=0$$
si $s>r$ et $n_s>0$. Posons $Y=\f{1}^{m_1}\cdots \f{r}^{m_r}$ et
$X=\e{1}^{n_1}\cdots\e{s}^{n_s-1}$. Nous avons
$$\varphi (Y\otimes X\e{s})=\sum_{(Y)}\varphi (Y_{(1)}\otimes\e{s})
\varphi (Y_{(2)}\otimes X)$$
d'apr\`es \eqref{accoup2}. Or, d'apr\`es les propositions \ref{iso} et 
\ref{delta(X)}, on a
$$\Delta (Y)=\f{r}\otimes y+1\otimes z+\sum\left(\prod_{i\leq r}
\f{i}\right)\otimes t$$
o\`u $y,z,t\in U_-$.
On en d\'eduit que $\varphi (Y_{(1)}\otimes\e{s})=0$ d'apr\`es les propositions
\ref{lem2_dx} et \ref{lem6}.\end{proof}

Pour $j=1,2,3$, posons $\tilde{H}_j=\frac{h}{2}\sum_{i=1}^3b_{ij}H'_i$ o\`u les 
scalaires  $b_{ij}$ sont d\'efinis par \eqref{CetB}.

\begin{lemm}\label{lem14}\label{lem15}
On a $\varphi (\tilde{H}_j\otimes H_i)=\delta_{ij}$ pour $1\leq i,j\leq 3$ et
$$\varphi \left(\frac{\tilde{H}_1^{q_1}}{q_1!}\frac{\tilde{H}_2^{q_2}}{q_2!}
\frac{\tilde{H}_3^{q_3}}{q_3!}\otimes H_1^{p_1}H_2^{p_2}H_3^{p_3}\right)
=\delta_{p_1q_1}\delta_{p_2q_2}\delta_{p_3q_3}.$$
\end{lemm}
\begin{proof} Laiss\'ee au lecteur.\end{proof}

\begin{prop}\label{duale}
Les familles
$$E^{n_1}_{\beta_1}\dots E^{n_7}_{\beta_7}H_1^{\ell_1}H_2^{\ell_2}H_3^{\ell_3}
\;\;\text{et} \;\;
(-1)^\xi\frac{F_{\beta_1}^{n_1}}
{(n_1)_1!\varphi_1^{n_1}}\cdots \frac{F_{\beta_7}^{n_7}}
{(n_7)_7!\varphi_7^{n_7}}\frac{\tilde{H}_1^{\ell_1}}
{\ell_1!}\frac{\tilde{H}_2^{\ell_2}}{\ell_2!}\frac{\tilde{H}_3^{\ell_3}}
{\ell_3!}$$
o\`u $\xi=n_6(n_2+n_3+n_5)+n_5(n_2+n_3)+n_3n_2$,
$\ell_i,\, n_j$ parcourent $\N$ pour $i=1,2,3$, 
$j=1,4,7$, et o\`u $n_i$ d\'ecrit $\{ 0,1\}$ pour
$i=2,3,5,6$, forment respectivement des $\C ((h))$-bases de 
$U_+\otimes_{\scriptscriptstyle \C [[h]]}\C ((h))$ et 
$U_-\otimes_{\scriptscriptstyle \C [[h]]}\C ((h))$ duales pour $\varphi$.
\end{prop}
\begin{proof} Soient $n_i,m_i,p_j,q_j\in\N$ ($i=1,4,7$, $j=1,2,3$) et 
$n_k,m_k\in\{0,1\}$ ($k=2,3,5,6$). D'apr\`es les lemmes \ref{lem1}, 
\ref{phi(F_E)} et \ref{lem14}, on a 
\begin{multline*}
\varphi \left(\!\frac{F_{\beta_1}^{m_1}}
{(m_1)_1!\varphi_1^{m_1}}\cdots \frac{F_{\beta_7}^{m_7}}
{(m_7)_7!\varphi_7^{m_7}}\frac{\tilde{H}_1^{q_1}}
{q_1!}\frac{\tilde{H}_2^{q_2}}{q_2!}\frac{\tilde{H}_3^{q_3}}{q_3!}\otimes
E_{\beta_1}^{n_1}\cdots E_{\beta_7}^{n_7}H_1^{p_1}H_2^{p_2}H_3^{p_3}\right) \\
=(-1)^\xi\prod_{\overset{1\leq i\leq 7}{\scriptscriptstyle 1\leq j\leq 3}}
\delta_{n_im_i}\delta_{p_jq_j},
\end{multline*}
o\`u $\xi=n_6(n_2+n_3+n_5)+n_5(n_2+n_3)+n_3n_2$. La proposition d\'ecoule
alors des propositions~\ref{iso} et~\ref{commut}.\end{proof}

\begin{proof}[D\'emonstration du th\'eor\`eme \ref{theoR}] D'apr\`es la
proposition~\ref{duale}, l'\'el\'ement
\begin{multline*}
R_{{\mathcal D}}=\sum_{\overset{m_1,m_2,m_3,n_1,n_4,n_7\in\N}{\scriptscriptstyle 
n_2,n_3,n_5,n_6\in\{0,1\}}}
\frac{(-1)^{n_6(n_5+n_3+n_2)+n_5(n_3+n_2)+n_3n_2}}{(n_1)_1!(n_4)_4!(n_7)_7!
m_1!m_2!m_3!\varphi_1^{n_1}\cdots\varphi_7^{n_7}}\\
\times\e{1}^{n_1}\cdots\e{7}^{n_7}
H_1^{m_1}\cdots H_3^{m_3}\otimes\f{1}^{n_1}\cdots\f{7}^{n_7}\tilde{H}_1^{m_1}
\cdots\tilde{H}_3^{m_3}
\end{multline*}
est une $R$-matrice universelle pour ${\mathcal D}$ ({\em cf.} 
\cite{drinfeld85,ZGB,rosso89bis}). 
On en d\'eduit que l'image 
$R\in{\overline{\mathcal D}}\otimes{\overline{\mathcal D}}$ de 
$R_{{\mathcal D}}$ par la projection canonique de ${\mathcal D}$ sur ${\overline{\mathcal D}}$ est une $R$-matrice 
universelle pour ${\overline{\mathcal D}}$. Or,
\begin{multline*}
(-1)^{\scriptscriptstyle n_6(n_5+n_3+n_2)+n_5(n_3+n_2)+n_3n_2}\e{1}^{n_1}\cdots
\e{7}^{n_7}\otimes\f{1}^{n_1}\cdots\f{7}^{n_7} \\
=(\e{1}^{n_1}\otimes\f{1}^{n_1})\cdots (\e{7}^{n_7}\otimes\f{7}^{n_7}),
\end{multline*}
et
$$\sum_{m_1,m_2,m_3\in\N}\frac{H_1^{m_1}H_2^{m_2}H_3^{m_3}\otimes
\tilde{H}_1^{m_1}\tilde{H}_2^{m_2}\tilde{H}_3^{m_3}}{m_1!\, m_2!\, m_3!}=
\exp\left(\frac{h}{2}\Big(\sum_{i,j}b_{ij}H_i\otimes H_j\Big)\right),$$
donc 
$$R=\exp_1\left(\frac{E_{\beta_1}\otimes F_{\beta_1}}{\varphi_1}\right)\cdots 
\exp_7\left(\frac{E_{\beta_7}\otimes F_{\beta_7}}{\varphi_7}\right)
\exp\left(\frac{h}{2}\Big(\sum_{i,j}b_{ij}H_i\otimes H_j\Big)\right)$$
d'apr\`es les formules \eqref{ei_et_phii}. La proposition
\ref{isomorphes} permet achever la d\'emonstration.\end{proof}


\section{D\'emonstration du th\'eor\`eme~\ref{theo2}}\label{demo2}
Dans ce paragraphe, nous d\'emontrons le th\'eor\`eme~\ref{theo2} \'enonc\'e au
\S\ref{resultats}. Nous commençons par d\'efinir un supermodule de rang six,
puis nous rappelons la construction de la cat\'egorie ruban\'ee associ\'ee; nous
terminerons par la d\'emonstration du th\'eor\`eme.

Dans tout le paragraphe~\ref{demo2}, le scalaire $x$ qui entre dans la d\'efinition
de $\hdx$ vaut $x=1$.

\subsection{Un supermodule de rang six}\label{par:module}

Soit $M_0$ le $C[[h]]$-module topologiquement libre de base $(v_1,v_2)$
et $M_1$ le $C[[h]]$-module topologiquement libre de base $(v_3,v_4,v_5,v_6)$.
Nous consid\'erons le supermodule $M=M_0\oplus M_1$ o\`u $M_0$ est la partie
paire et $M_1$ la partie impaire. On pose pour $i=2,3$, $j=1,2$,
\begin{gather*}
H_2v_1=H_3v_1=v_1,\quad H_2v_2=H_3v_2=v_2,\quad
E_iv_j=F_iv_j=0,\\
H_1v_1=v_1,\quad H_1v_2=-v_2,\quad
E_1v_1=F_1v_2=0,\quad F_1v_1=v_3,\quad E_1v_2=-v_6,
\end{gather*}
\begin{gather*}
H_2v_3=H_3v_3=v_3,\quad H_2v_6=H_3v_6=-v_6,\\ H_2v_4=-H_3v_4=-v_4,\quad
H_2v_5=-H_3v_5=v_5,\\
E_2v_3=E_3v_3=E_2v_5=E_3v_4=F_2v_6=F_3v_6=F_2v_4=F_3v_5=0,\\
E_2v_4=E_3v_5=v_3,\quad E_2v_6=v_5,\quad E_3v_6=v_4,\\
F_2v_5=F_3v_4=v_6,\quad F_2v_3=v_4,\quad F_3v_3=v_5,\\
E_1v_4=E_1v_5=E_1v_6=F_1v_3=F_1v_4=F_1v_5=0,\\
E_1v_3=v_1,\quad F_1v_6=v_2,\quad H_1v_3=v_3,\quad H_1v_4=H_1v_5=0,\quad
H_1v_6=-v_6.
\end{gather*}

\begin{prop}\label{prop:module}
Les formules pr\'ec\'edentes munissent $M$ d'une
structure de $\hdx$-module simple topologiquement libre isomorphe \`a son dual.
\end{prop}
\begin{proof} Nous laissons le soin au lecteur de v\'erifier que ces relations 
d\'efinissent une structure de $\hdx$-module. Montrons que $M$ est
simple. La sous-alg\`ebre de Hopf $U'$ de $\hdx$ topologiquement
engendr\'ee par $H_i,E_i,F_i$ pour $i=2,3$, est  isomorphe \`a
$U_h(sl_2)\hat{\otimes} U_h(sl_2)$.
En tant que $U'$-module, 
on a $M_0\cong \C[[h]]\oplus\C[[h]]$ et $M_1\cong V_1\hat{\otimes} V_1$. Comme 
$V_1\hat{\otimes} V_1$ est un
$U_h(sl_2)\hat{\otimes} U_h(sl_2)$-module simple, il suffit de v\'erifier que chacun 
des vecteurs $v_1,v_2,v_3$ engendre le module $M$ tout entier. Or, ceci est 
clair
d'apr\`es la forme des actions de $E_1$ et $F_1$. De plus, le lecteur v\'erifiera que
l'application $\alpha :M\rightarrow M^*$ d\'efinie par 
\begin{gather*}
\alpha (v_1)=-q^{-3}v^2,\;\; \alpha (v_2)=q^{-1}v^1,\;\; 
\alpha (v_3)=q^{-2}v^6,\\
\alpha (v_4)=-q^{-1}v^5,\;\; \alpha (v_5)=-q^{-1}v^4,\;\; \alpha (v_6)=v^3,
\end{gather*}
o\`u $(v^1,\dots ,v^6)$ d\'esigne la base duale de $(v_1,\dots ,v_6)$,
est un isomorphisme de modules.\end{proof}
Nous donnons maintenant le tressage  $c_{\scriptscriptstyle M,M}$ de $M$
induit par la $R$-matrice universelle de $\hdx$ via la 
formule~\eqref{twist_et_tressage}. 
L'automorphisme $c_{\scriptscriptstyle M,M}$   laisse  stable les parties 
paire et impaire de
$M^{\hat{\otimes} 2}$ car $R$ est un \'el\'ement  pair de $\hdx^{\hat{\otimes} 2}$. 
Le calcul de $c_{\sst M,M}$ a \'et\'e fait \`a l'aide de Maple. Nous
donnons les matrices de $c_0=c_{\scriptscriptstyle M,M|(M\hat{\otimes} M)_0}$ (resp.
$c_1=c_{\scriptscriptstyle M,M|(M\hat{\otimes} M)_1}$) dans la base $(v_i\otimes
v_j)_{\scriptscriptstyle \overset{1\leq i,j\leq 2}{\text{ou } 3\leq i,j\leq
6}}$ (resp.$(v_i\otimes v_j)_{\scriptscriptstyle \overset{1\leq i \leq 2,\,
3\leq j\leq 6}{\text{ou } 3\leq i\leq 6,\, 1\leq j\leq 2}}$) ordonn\'ee
suivant l'ordre lexicographique. Dans ces matrices, on pose $\lambda=q-q^{-1}$.
On v\'erifie que les polyn\^omes caract\'eristique
et minimal de $c_{\scriptscriptstyle M,M}$ sont
$(X-q)^{17}(X+q)(X+q^{-1})^{18}$ et $(X+q)(X-q)(X+q^{-1})$, respectivement. 

\begin{figure}
$$\setcounter{MaxMatrixCols}{20}
c_0=\psmallmatrixp{
q\phantom{q}& 0& 0& 0& 0& 0& 0& 0& 0& 0& 0& 0& 0& 0& 0& 0& 0& 0& 0& 0\cr
0& 0& \frac{1}{q}& 0& 0& 0& 0& 0& 0& 0& 0& 0& 0& 0& 0& 0& 0& 0& 0& 0\cr 0&
\frac{1}{q}& q^3-\frac{1}{q}& 0& 0& 0& 0& q\lambda & 0& 0& -q^2\lambda & 0& 0& -q^2\lambda& 0& 0&
q^3\lambda& 0& 0& 0\cr 0& 0& 0& q& 0& 0& 0& 0& 0& 0& 0& 0& 0& 0& 0& 0& 0&
0& 0& 0\cr 0& 0& 0& 0& -\frac{1}{q}& 0& 0& 0& 0& 0& 0& 0& 0& 0& 0& 0& 0& 0& 0&
0\cr 0& 0& 0& 0& 0& 0& 0& 0& -1& 0& 0& 0& 0& 0& 0& 0& 0& 0& 0& 0\cr 0&
0& 0& 0& 0& 0& 0& 0& 0& 0& 0& 0& -1& 0& 0& 0& 0& 0& 0& 0\cr 0& 0&
-\frac{\lambda}{q} & 0& 0& 0& 0& 0& 0& 0& 0& 0& 0& 0& 0& 0& -q& 0& 0& 0\cr 0&
0& 0& 0& 0& -1& 0& 0& \lambda & 0& 0& 0& 0& 0& 0& 0& 0& 0& 0& 0\cr 0& 0&
0& 0& 0& 0& 0& 0& 0& -\frac{1}{q}& 0& 0& 0& 0& 0& 0& 0& 0& 0& 0\cr 0& 0&
\lambda & 0& 0& 0& 0& 0& 0& 0& 0& 0& 0& -q& 0& 0& q\lambda & 0& 0& 0\cr 0& 0&
0& 0& 0& 0& 0& 0& 0& 0& 0& 0& 0& 0& 0& 0& 0& -1& 0& 0\cr 0& 0& 0& 0&
0& 0& -1& 0& 0& 0& 0& 0& \lambda & 0& 0& 0& 0& 0& 0& 0\cr 0& 0& \lambda & 0&
0& 0& 0& 0& 0& 0& -q& 0& 0& 0& 0& 0& q\lambda & 0& 0& 0\cr 0& 0& 0& 0& 0&
0& 0& 0& 0& 0& 0& 0& 0& 0& -\frac{1}{q}& 0& 0& 0& 0& 0\cr 0& 0& 0& 0& 0& 0& 0&
0& 0& 0& 0& 0& 0& 0& 0& 0& 0& 0& -1& 0\cr 0& 0& -q\lambda & 0& 0& 0& 0&
-q& 0& 0& q\lambda & 0& 0& q\lambda & 0& 0& -q\lambda^2& 0& 0& 0\cr 0& 0& 0&
0& 0& 0& 0& 0& 0& 0& 0& -1& 0& 0& 0& 0& 0& \lambda & 0& 0\cr 0& 0& 0& 0&
0& 0& 0& 0& 0& 0& 0& 0& 0& 0& 0& -1& 0& 0& \lambda & 0\cr 0& 0& 0& 0& 0&
0& 0& 0& 0& 0& 0& 0& 0& 0& 0& 0& 0& 0& 0& -\frac{1}{q}
}$$
\end{figure}

\begin{figure}
$$\setcounter{MaxMatrixCols}{16}
c_1=\psmallmatrixg{
0& 0& 0& 0& 0& 0& 0& 0& 1& 0& 0& 0& 0& 0& 0& 0\cr 0& 0& 0&
0& 0& 0& 0& 0& 0& 0& 1& 0& 0& 0& 0& 0\cr 0& 0& 0& 0& 0& 0& 0& 0& 0& 0&
0& 0& 1& 0& 0& 0\cr 0& 0& 0& 0& 0& 0& 0& 0& 0& 0& 0& 0& 0& 0& 1& 0\cr
0& 0& 0& 0& \lambda & 0& 0& 0& 0& 1& 0& 0& 0& 0& 0& 0\cr 0& 0& 0& 0& 0&
\lambda & 0& 0& 0& 0& 0& 1& 0& 0& 0& 0\cr 0& 0& 0& 0& 0& 0& \lambda & 0& 0&
0& 0& 0& 0& 1& 0& 0\cr 0& 0& 0& 0& 0& 0& 0& \lambda & 0& 0& 0& 0& 0& 0&
0& 1\cr 1& 0& 0& 0& 0& 0& 0& 0& \lambda & 0& 0& 0& 0& 0& 0& 0\cr 0& 0& 0&
0& 1& 0& 0& 0& 0& 0& 0& 0& 0& 0& 0& 0\cr 0& 1& 0& 0& 0& 0& 0& 0& 0& 0&
\lambda & 0& 0& 0& 0& 0\cr 0& 0& 0& 0& 0& 1& 0& 0& 0& 0& 0& 0& 0& 0& 0&
0\cr 0& 0& 1& 0& 0& 0& 0& 0& 0& 0& 0& 0& \lambda & 0& 0& 0\cr 0& 0& 0& 0&
0& 0& 1& 0& 0& 0& 0& 0& 0& 0& 0& 0\cr 0& 0& 0& 1& 0& 0& 0& 0& 0& 0& 0&
0& 0& 0& \lambda & 0\cr 0& 0& 0& 0& 0& 0& 0& 1& 0& 0& 0& 0& 0& 0& 0&
0
}$$
\end{figure}

\subsection{Cat\'egorie ruban\'ee associ\'ee \`a $M$}\label{catrub}

Soit ${\mathcal C}_M$ la sous-cat\'egorie pleine de la cat\'egorie des $\hdx$-modules
dont les objets sont les produits tensoriels finis it\'er\'es de $M$ et $M^*$. 
Cette cat\'egorie
est une cat\'egorie tensorielle d'objet unit\'e $\C[[h]]$ et est munie d'une
dualit\'e  induite par les applications
$$b_{\sst M}:\C[[h]]\longrightarrow M\otimes M^*,\quad d_{\sst M}:M^*\otimes M
\longrightarrow \C[[h]]$$
d\'efinies par $b_{\sst M}(1)=\sum_{i=1}^6 v_i\otimes v^i$ et $d_{\sst M} 
(f\otimes
x)=f(x)$ avec  $x\in M$ et $f\in M^*$. Cette cat\'egorie est \'egalement ruban\'ee
({\em cf.} \S\ref{resultats}). Le carr\'e de l'inverse du twist est donn\'e par la
formule  (lemme XIV.3.4 de \cite{K}) 
$$\theta_{\sst M}^{-2}=(d_{\sst M}c_{\sst M,M^*}\otimes\id_{\sst
M})(\id_{\sst M}\otimes c_{\sst M,M^*}b_{\sst M}).$$
Comme $\alpha:M\rightarrow M^*$ est un isomorphisme, la naturalit\'e du tressage
donne 
$$c_{\sst M,M^*}=(\alpha^{-1}\otimes\id_{\sst M})c_{\sst M,M}(\id_{\sst
M}\otimes\alpha ).$$ 
Un calcul \`a l'aide de Maple montre que $\theta_{\sst
M}^2=q^{-2}\id_{\sst M}$. Par cons\'equent,
$$\theta_{\sst M}=q^{-1}\id_{\sst M},$$
car $\theta_{\sst M}\equiv\id_{\sst M} \mod h$. D\'efinissons les morphismes
\begin{gather*}
b'_{\sst M}=(\id_{\sst M^*}\otimes\theta_{\sst M})c_{\sst M,M^*}b_{\sst M}
:\C[[h]]\longrightarrow M^*\otimes M, \\
d'_{\sst M}=d_{\sst M}c_{\sst M,M^*}(\theta_{\sst M}\otimes\id_{\sst M^*})
:M\otimes M^*\longrightarrow \C[[h]],
\end{gather*}
et posons
\begin{equation}
\begin{gathered}
d=d_{\sst M}(\alpha\otimes\id_{\sst M})
:M\otimes M\longrightarrow \C[[h]],\\
b=(\id_{\sst M}\otimes\alpha^{-1})b_{\sst M}
:\C[[h]]\longrightarrow M\otimes M, \\
d'=d'_{\sst M}(\id_{\sst M}\otimes\alpha),\qquad 
b=(\alpha^{-1}\otimes\id_{\sst M})b_{\sst M}.
\end{gathered}
\end{equation}

\begin{lemm}\label{lem:egalites}
Nous avons les \'egalit\'es suivantes :
\begin{gather*}
c_{\sst M,M}-c_{\sst M,M}^{-1}=(q-q^{-1})(\id_{\sst M\otimes M}-bd), \\
db(1)=2,\qquad (\id_{\sst M}\otimes d)(c_{\sst M,M}\otimes\id_{\sst
M})(\id_{\sst M}\otimes b)=-q^{-1}\id_{\sst M},\\
b'=-b,\qquad d'=-d.
\end{gather*}
\end{lemm}
\begin{proof} Ces relations peuvent se v\'erifier \`a l'aide de Maple. On peut 
\'egalement montrer que les espaces propres de $c_{\sst M,M}$ sont des 
sous-modules simples de $M\otimes M$.\end{proof}

\subsection{Invariant d'entrelacs associ\'e}\label{invariant}

Nous donnons dans ce paragraphe la d\'emonstration du th\'eor\`eme~\ref{theo2}.
Notons ${\mathcal E}$ la cat\'egorie dont les objets sont les suites finies de
$+$ et de $-$, y compris la suite vide, et dont les morphismes sont engendr\'es
par les enchev\^etrements orient\'es. Pour plus de d\'etails sur cette cat\'egorie,
on pourra consulter~\cite{Tu2}.

La cat\'egorie ${\mathcal C}_{M}$ \'etant ruban\'ee, on sait qu'il existe un
unique foncteur ${\mathcal F}:{\mathcal E}\longrightarrow {\mathcal C}_M$ tel que
\begin{gather*}
{\mathcal F}(+)=M,\qquad {\mathcal F}(-)=M^*, \\
{\mathcal F}(\dessin{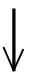})=\id_{\sst M},\qquad
{\mathcal F}(\dessin{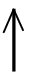})=\id_{\sst M^*},\\
{\mathcal F}(\dessin{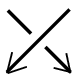})=c_{\sst M,M}, \quad 
{\mathcal F}(\dessin{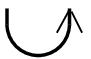})=b_{\sst M},\quad {\mathcal
F}(\dessin{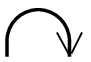})=d_{\sst M},
 \end{gather*}
et que si $L$ est un entrelacs orient\'e parall\'elis\'e, alors 
$${\mathcal F}(L)\in \End_{\sst\hdx}(\C[[h]])\cong\C[[h]]$$ 
d\'efinit un invariant d'entrelacs.

Avant de d\'emontrer le th\'eor\`eme~\ref{theo2}, nous avons besoin des deux
lemmes suivants. Soit $L$ est un entrelacs orient\'e parall\'elis\'e. Fixons une
repr\'esentation planaire de $L$ et soit $n_+$ (resp. $n_-$) la somme du nombre
de $\dessin{cup+}$ et de $\dessin{cap+}$ (resp. $\dessin{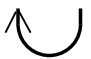}$ et de
$\dessin{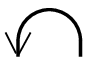}$) appara\^issant dans cette repr\'esentation.

\begin{lemm}\label{le_signe}
La parit\'e des entiers ${n_+}$  et ${n_-}$ ne d\'epende que de la classe
d'isotopie de $L$. De plus, nous avons ${n_+}\equiv {n_-}\mod 2$ .
\end{lemm}
\begin{proof} La premi\`ere affirmation d\'ecoule du fait que les mouvements de
Reidemeisters (pour les entrelacs fram\'es) ne changent pas la parit\'e de  
$n_+$ et $n_-$. Il en est de m\^eme de l'\'egalit\'e
$$\dessin{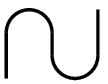}=\dessin{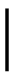},$$
(orientation quelconque).

Pour d\'emontrer l'\'egalit\'e des congruences, consid\'erons comme repr\'esentation
planaire de $L$ la cl\^oture d'une tresse, \`a laquelle nous ajoutons, suivant
la valeur du framing de $L$, un certain nombre de boucles.  Il est clair que
pour cette repr\'esentation, on a $n_+=n_-$.\end{proof}

Nous posons $\varepsilon(L)=(-1)^{n_+}$.
\begin{lemm}\label{pasorient}
L'invariant ${\mathcal F}(L)$ est ind\'ependant de l'orientation de $L$.
\end{lemm}
\begin{proof} Notons $\tilde{L}$ l'entrelacs $L$ dont l'orientation d'une de ses
composantes est invers\'ee.
Alors, d'apr\`es les
lemmes~\ref{lem:egalites} et~\ref{le_signe}, on a
$${\mathcal F}(L)=\varepsilon(L)\tilde{\mathcal F}(L),$$
o\`u $\tilde{\mathcal F}(L)\in\End_{\sst\hdx}(\C[[h]])$ est le morphisme obtenu
\`a partir de ${\mathcal F}(L)$ en rempla\c{c}ant $b_{\sst M}$ et $b'_{\sst
M}$ (resp. $d_{\sst M}$ et $d'_{\sst M}$)  par $b$ (resp. $d$), ainsi que
$c_{\sst M,M^*}$, $c_{\sst M^*,M}$, $c_{\sst M^*,M^*}$ (resp. $c^{-1}_{\sst
M,M^*}$,  $c^{-1}_{\sst M^*,M}$, $c^{-1}_{\sst M^*,M^*}$) par
$c_{\sst M,M}$ (resp. $c^{-1}_{\sst M,M}$). Comme $\varepsilon(L)=
\varepsilon(\tilde{L})$ et $\tilde{\mathcal F}(L)=\tilde{\mathcal F}(\tilde{L})$, le
lemme est d\'emontr\'e.\end{proof}

\begin{proof}[D\'emonstration du th\'eor\`eme~\ref{theo2}] Nous reprenons les notations
d\'efinies dans la d\'emonstration du lemme pr\'ec\'edent, et nous posons
$${\mathcal I}_L=\varepsilon(L){\mathcal F}(L)=\tilde{\mathcal F}(L).$$
Les lemmes~\ref{le_signe} et~\ref{pasorient} assurent que ${\mathcal I}_L$ est
bien un invariant d'un entrelacs parall\'elis\'e non orient\'e. De plus, 
d'apr\`es le lemme~\ref{lem:egalites}, cet
invariant v\'erifie bien les relations skeins du polyn\^ome de Dubrovnik,
dans lequel on a pos\'e $a=-q^{-1}$ et $z=q-q^{-1}$.
\end{proof}


\backmatter
\bibliography{dx}
\bibliographystyle{smfplain}

\end{document}